\documentclass[a4paper, 11pt, twoside, final]{article}

\usepackage{enumerate,graphicx,amsfonts,amsmath,psfrag}
\usepackage{a4wide}

\numberwithin{equation}{section}
\allowdisplaybreaks[4]

\newtheorem{proposition}{Proposition}[section]
\newtheorem{theorem}[proposition]{Theorem}
\newtheorem{lemma}[proposition]{Lemma}
\newtheorem{corollary}[proposition]{Corollary}
\newtheorem{definition}[proposition]{Definition}

\newenvironment{proof}{\smallskip\noindent\emph{\textbf{Proof.}}\hspace{1pt}}%
{\hspace{-5pt}{\nobreak\quad\nobreak\hfill\nobreak$\square$\vspace{8pt}%
\par}\smallskip\goodbreak}

\newenvironment{proofof}[1]{\smallskip\noindent\emph{\textbf{Proof of #1.}}%
\hspace{1pt}}{\hspace{-5pt}{\nobreak\quad\nobreak\hfill\nobreak%
$\square$\vspace{8pt}\par}\smallskip\goodbreak}

\renewcommand{\div}{\mathinner{\mathrm{div}}}

\renewcommand{\L}[1]{\mathbf{L^#1}}

\newcommand{\Wloc}[2]{\mathbf{W^{#1,#2}_{loc}}}
\newcommand{\C}[1]{\mathbf{C^{#1}}}
\newcommand{\Cc}[1]{\mathbf{C_c^{#1}}}

\newcommand{\W}[2]{\mathbf{W^{#1,#2}}}
\newcommand{\modulo}[1]{{\left|#1\right|}}
\newcommand{\norma}[1]{{\left\|#1\right\|}}
\newcommand{\reali}{{\mathbb{R}}}
\newcommand{\naturali}{{\mathbb{N}}}
\newcommand{\tv}{\mathrm{TV}}
\newcommand{\BV}{\mathbf{BV}}

\renewcommand{\epsilon}{\varepsilon}
\renewcommand{\phi}{\varphi}

\renewcommand{\d}[1]{\mathinner{\mathrm{d}{#1}}}

\newcommand{\Caption}[1]{\caption{\small#1}}
\begin{document}

\title{An Analytical Framework to Describe the\\Interactions Between
  Individuals and a Continuum}

\author{Rinaldo M.~Colombo$^1$ \and Magali L\'ecureux-Mercier$^2$}

\footnotetext[1]{Department of Mathematics, Brescia University, Via
  Branze 38, 25133 Brescia, Italy}

\footnotetext[2]{Universit\'e d'Orl\'eans, UFR Sciences, B\^atiment de
  math\'ematiques - Rue de Chartres B.P. 6759 - 45067 Orl\'eans cedex
  2, France}

\maketitle

\begin{abstract}
  \noindent We consider a discrete set of individual agents
  interacting with a continuum. Examples might be a predator facing a
  huge group of preys, or a few shepherd dogs driving a herd of
  sheeps. Analytically, these situations can be described through a
  system of ordinary differential equations coupled with a scalar
  conservation law in several space dimensions. This paper provides a
  complete well posedness theory for the resulting Cauchy problem. A
  few applications are considered in detail and numerical integrations
  are provided.

  \medskip


  \medskip


\end{abstract}

\section{Introduction}

In various situations a small set of individuals has to interact with
a continuum. A first famous examples comes from the fairy tale of the
pied piper~\cite{GrimmBrothers}, where a musician frees a city from
rats using his magic flute. An entirely different case is that of
shepherd dogs confining sheeps while pasturing, or that of a wild
predator seeking to split a flock of preys. From a deterministic point
of view, studying these phenomena leads to a dynamical system
consisting of ordinary differential equations for the evolution of the
agents and partial differential equations for that of the
continuum. Here, motivated by the present toy applications, we choose
scalar conservation laws for the description of the continuum's
evolution. In particular, no diffusion is here considered. On one
side, this choice makes the analytical treatment technically more
difficult, due to the possible singularities arising in the density
that describes the continuum. On the other hand, we obtain a framework
where all propagation speeds are finite. As a consequence, for
instance, a continuum initially confined in a bounded region will
remain in a (larger but) bounded region at any positive time. This
allows to state problems concerning the support of the continuum, such
as confinement problems (the rats should leave the city, or the
shepherd dogs should keep sheeps inside a given area) or far more
complex ones (how can a predator split the support of the density of
its preys?).

In the current literature, similar problems have been considered with
a great variety of analytical tools, see for
instance~\cite{BressanDeLellis} for a fire confinement problem modeled
through differential inclusions, or~\cite{CapassoMichelettiMorale} for
a tumor--induced angiogenesis described through a stochastic geometric
model. Other examples are provided by the interaction of a fluid
(liquid or gas) with a rigid body or with an elastic structure, like a
membrane, see~\cite{Serre1987, VasquezZuazua}: the evolution of the
rigid body is described by a system of ordinary differential
equations, while the evolution of the fluid is subject to partial
differential equations like Navier-Stokes or Euler equations. Further
results are currently available in the 1D case. For instance, a
problem motivated by traffic flow is considered
in~\cite{LattanzioMauriziPiccoli}; the piston problem, a blood
circulation model and a supply chain model are considered
in~\cite{BorscheColomboGaravello}.

Formally, we are thus lead to the dynamical system
\begin{equation}
  \label{eq:Problem}
  \left\{
    \begin{array}{l}
      \partial_t \rho
      +
      \div_x f \left(t, x, \rho, p(t)\right)
      =
      0
      \\
      \dot p
      =
      \phi\left( t, p, \left(A \rho(t)\right) (p) \right)
      \\
      \rho(0,x) = \bar \rho(x)
      \\
      p(0) = \bar p
    \end{array}
  \right.
  \qquad
  \begin{array}{rcl}
    (t, x )& \in & \reali^+ \times \reali^{N_x}
    \\
    \rho & \in & \reali^+
    \\
    p & \in & \reali^{N_p}
  \end{array}
\end{equation}
where the unknowns are $\rho$ and $p$. The former one, $\rho =
\rho(t,x)$ is the density describing the macroscopic state of the
continuum while the latter, $p = p(t)$, characterizes the state of the
individuals. It can be for instance the vector of the individuals'
positions or of the individuals' positions and speeds. The dynamics of
the continuum is described by the flow $f$, which in general can be
thought as the product $f = \rho \, v$ of the density $\rho$ and a
suitable speed $v = v(t, x, \rho, p)$. The vector field $\phi$ defines
the dynamics of the individuals at time $t$ and it depends from the
continuum density $\rho(t)$ through a suitable average $A
\left(\rho(t)\right)$. Our driving example below is the convolution in
the space variable $A \left(\rho(t) \right) = \rho(t) * \eta$, with a
smooth compactly supported kernel $\eta$.

Below we address and solve the first mathematical questions that arise
about~\eqref{eq:Problem}, i.e.~the existence and uniqueness of entropy
solutions, their stability with respect the data and the equation, and
the existence of optimal controls. A first well posedness result, that
applies to general initial data, is provided in
Theorem~\ref{thm:main}. As usual in this context, see
also~\cite{ColomboHertyMercier, ColomboMercierRosini, Kruzkov,
  MercierStability}, the hypotheses on $f$ are rather
intricate. However, the present framework naturally applies to
situations in which the continuum can be supposed initially confined
in a bounded region, i.e.~$\rho$ vanishes outside a compact subset of
$\reali^{N_x}$. In this case, Corollary~\ref{cor:compact} applies and
the hypotheses on $f$ are greatly simplified.

\smallskip

The next section presents the analytical well-posedness
results. Section~\ref{sec:Appl} is devoted to various applications,
while all proofs are deferred to the last section.

\section{Notation and Analytical Results}

We now collect the various assumptions on~(\ref{eq:Problem}) that
allow us to prove well posedness, i.e.~the existence of solutions,
their uniqueness and their stability with respect to data and
equations. The hypotheses collected below are essentially those that
ensure the well posedness of the conservation law and, separately, of
the ordinary differential equation.

Throughout, we denote $\reali^+ = \left[0, +\infty \right[$ and
$B_{\reali^{N_p}}(x,r)$ denotes the closed ball in $\reali^{N_p}$
centered at $x$ with radius $r$.  Let $T_{\max} \in [0, +\infty]$ and
call $I=[0, T_{\max}]$ if $T_{\max} < +\infty$, while $I=\reali^+$
otherwise.  The real parameter $R$, i.e.~the maximal possible density
is fixed and positive. For a given compact set $K$ in $\reali^{N_p}$
and a $T>0$, we denote $\Omega_T = [0,T] \times \reali^{N_x} \times
[0,R] \times K$.

\medskip

\noindent\textbf{Flow of the continuum:} at point $x$ and time $t$,
the continuum flows with a flux $f = f\left(t, x, \rho(t,x), p(t)
\right)$ that depends on time $t$, on the space variable $x$, on the
continuum density $\rho$ evaluated at $(t,x)$ and on the state $p$ of
the individuals at time $t$. We require the following regularity:
\begin{description}
\item[(f)] The flow $f \colon I \times \reali^{N_x} \times [0,R]
  \times \reali^{N_p} \to \reali^{N_x}$ is such that
  \begin{enumerate}[\bf(f.1)]
  \item \label{it:f1} $f \in \C2 (I \times \reali^{N_x} \times [0,R]
    \times \reali^{N_p}; \reali^{N_x})$.
  \item \label{it:f2} For all $(t, x, p) \in I \times \reali^{N_x}
    \times \reali^{N_p}$, $f(t, x, 0, p) = f(t, x, R, p) = 0$.
  \item \label{it:f3} For all $T \in I$ and for all compact subsets $K
    \subset \reali^{N_p}$, there exists a constant $C_f$ such that for
    $t \in [0,T]$, $x \in \reali^{N_x}$, $\rho \in [0,R]$ and $p \in
    K$,
    \begin{displaymath}
      \norma{\partial_\rho f(t, x, \rho, p)} < C_f\,,
      \modulo{\div_x f(t, x, \rho, p)} < C_f\,.
    \end{displaymath}
    % \item \label{it:f4} For all $T \in I$ and for all compact
    %   subsets
    %   $K \subset \reali^{N_p}$, there exists a constant $C$ such
    %   that
    %   for
    %   $t \in [0,T]$, $x \in \reali^{N_x}$, $\rho \in [0,R]$ and $p
    %   \in
    %   K$,
    %   $\modulo{\div_x f(t, x, \rho, p)} < C$.
  \item \label{it:f5} For all $T \in I$ and for all compact subsets $K
    \subset \reali^{N_p}$, there exists a constant $C_f$ such that for
    $t \in [0,T]$, $x \in \reali^{N_x}$, $\rho \in [0,R]$ and $p \in
    K$,
    \begin{displaymath}
      \norma{\nabla_x \partial_\rho f(t, x, \rho, p)} <  C_f\,.
    \end{displaymath}
  \item \label{it:f6} For all compact subsets $K \subset
    \reali^{N_p}$, there exists a constant $C_f$ such that\break
    \begin{displaymath}
      \int_I \int_{\reali^{N_x}} \sup_{p \in K, \rho \in [0,R]}
      \norma{\nabla_x \div_x f(t, x, \rho, p)} \d{x} \d{t} < C_f\,,
    \end{displaymath}
  \item \label{it:f7} For all compact subsets $K \subset
    \reali^{N_p}$, there exists a constant $C_f$ such that\break
    \begin{displaymath}
      \int_I \int_{\reali^{N_x}} \sup_{p \in K, \rho \in [0,R]}
      \norma{\div_x f(t, x, \rho, p)} \d{x} \d{t} < C_f\,.
    \end{displaymath}
  \item \label{it:f8} For all $T \in I$ and for all compact subsets $K
    \subset \reali^{N_p}$, there exists a constant $C_f$ such that for
    $t \in [0,T]$, $\rho \in [0,R]$ and $p \in K$,
    \begin{displaymath}
      \int_{\reali^{N_x}}
      \norma{\nabla_p\div_x f(t, x, \rho, p)}
      \d{x}
      <
      C_f\,,
      \norma{\nabla_p\partial_\rho f(t, x, \rho, p)}
      <
      C_f \mbox{ for all }x \in \reali^{N_x}\,.
    \end{displaymath}
  \end{enumerate}
\end{description}
\noindent Condition~\textbf{(f.\ref{it:f2})} states that at the
maximal density $\rho = R$, the continuum is at congestion and can not
move. Assumption~\textbf{(f.\ref{it:f2})} has a key importance. The
first part ensures the finite propagation speed of the solution to the
partial differential equation, see~Proposition~\ref{prop:FiniteSpeed}
or~\cite[Theorem~1]{Kruzkov}. The second part ensures that the
solutions are bounded, similarly to the role of the
sublinearity~\textbf{($\boldsymbol{\phi}$.3)} in the ordinary
differential equation.

All these assumptions are satisfied, for instance, by vector fields of
the form $u(\rho, x, p) = v(\rho) \, \vec{\textbf{v}}(x,p)$ with $v
\in \C2([0, R]; \reali)$ and $\vec{\textbf{v}} \in \Cc2(\reali^{N_x}
\times \reali^{N_x}; \reali^{N_x})$.

We note that if $f$ does not depend explicitly on $t$ and $x$, which
is a usual situation when dealing with systems of conservation laws in
one space dimension, then the above assumptions reduce to
only~\textbf{(f.\ref{it:f1})}, \textbf{(f.\ref{it:f2})}, the first
part of~\textbf{(f.\ref{it:f3})} and the first part
of~\textbf{(f.\ref{it:f8})}.

Moreover, Corollary~\ref{cor:compact} shows that whenever the initial
density distribution $\bar\rho$ has compact support, then the
requirements on $f$ are reduced, since only~\textbf{(f.\ref{it:f1})},
\textbf{(f.\ref{it:f2})} and~\textbf{(f.\ref{it:f3})} are necessary.

\medskip

\noindent\textbf{Speed of the individuals:} at time $t$,
the individuals' state changes with a speed $\phi = \phi\left(t, p(t),
  A\left(\rho(t)\right) \left(p(t)\right)\right)$ that depends on time
$t$, on the individuals' state $p$ at time $t$ and on an average
$A\left(\rho(t)\right)$ of the continuum density $\rho$ evaluated at
time $t$ and computed at $p(t)$. On the averaging operator $A$ we
require the following conditions.
\begin{description}
\item [(A)] $A \colon \L1(\reali^{N_x}; \reali) \to
  \W{1}{\infty}(\reali^{N_p};\reali^{N_r})$ is linear and continuous,
  i.e.~there exists a constant $C_A$ such that for all $\rho \in
  \L1(\reali^{N_x};\reali)$
  \begin{displaymath}
    \norma{A \rho}_{\W{1}{\infty}}
    \leq
    C_A \; \norma{\rho}_{\L1} \,.
  \end{displaymath}
\end{description}
\noindent Below, the operator norm of $A$ is denoted
$\norma{A}_{\mathcal{L}(\L1,\W{1}{\infty})}$. For instance, in the
case $N_p = N_x$, a typical example of such an operator $A$ is
$\left(A(\rho)\right)(p) = (\rho * \eta) (p)$ for a kernel $\eta \in
\Cc1(\reali^{N_x};\reali)$ with $\int_{\reali^{N_p}} \eta \d{x}=1$.

The speed law $\phi$ satisfies the assumptions:
\begin{description}
\item[($\boldsymbol{\phi}$)] The vector field $\phi \colon \reali^+
  \times \reali^{N_p} \times \reali^{N_r} \to \reali^{N_p}$ is such
  that
  \begin{enumerate}[\bf($\boldsymbol{\phi}$.1)]
  \item $t \mapsto \phi(t,p,r)$ is measurable for all $p \in
    \reali^{N_p}$ and all $r \in \reali^{N_r}$;
  \item there exists a function $C_\phi \in \L1 (I;\reali^+)$ such
    that for a.e.~$t \in I$, $p_1,p_2 \in \reali^{N_p}$ and $r_1,r_2
    \in \reali^{N_r}$,
    \begin{displaymath}
      \norma{\phi(t,p_1,r_1) - \phi(t, p_2, r_2)}
      \leq
      C_\phi(t)
      \left(
        \norma{p_1-p_2} + \norma{r_1-r_2}
      \right) \,;
    \end{displaymath}
  \item there exists a function $C_\phi \in \L1 (I;\reali^+)$ such
    that for a.e.~$t \in [0,T]$, for all $p \in \reali^{N_p}$ and for
    all $r \in \reali^{N_r}$,
    \begin{displaymath}
      \norma{\phi(t,p,r)}
      \leq
      C_\phi(t) \left( 1 + \norma{p} \right)\,.
    \end{displaymath}
  \end{enumerate}
\end{description}
\medskip

\noindent These hypotheses are motivated by the standard theory of
Caratheodory ordinary differential equations,
see~\cite[\S~1]{Filippov}. All the above assumptions~\textbf{(f)},
\textbf{(A)} and~\textbf{($\boldsymbol{\phi}$)} are satisfied in the
applications considered in Section~\ref{sec:Appl}.

As a first step in the analytical treatment of~(\ref{eq:Problem}), we
rigorously state what we mean by \emph{solution}
to~(\ref{eq:Problem}).

\begin{definition}
  \label{def:sol}
  Fix $\bar \rho \in (\L1\cap \BV) \left( \reali^{N_x}; [0,R] \right)$
  and $\bar p \in \reali^{N_p}$. A pair $(\rho, p)$ with
  \begin{displaymath}
    \rho \in
    \C0 \left(I; \L1 ( \reali^{N_x}; [0,R]) \right)
    \quad \mbox{ and } \quad
    p \in
    \W{1}{1} (I; \reali^{N_p})
  \end{displaymath}
  is a solution to~(\ref{eq:Problem}) with initial datum $(\bar \rho,
  \bar p)$ if
  \begin{enumerate}[(i)]
  \item the map $\rho = \rho(t,x)$ is a Kru\v zkov solution to the
    scalar conservation law
    \begin{equation}
      \label{eq:SolHCL}
      \partial_t \rho
      +
      \div_x f\left( t, x, \rho, p(t) \right)
      =0
    \end{equation}
  \item the map $p = p(t)$ is a Caratheodory solution to the ordinary
    differential equation
    \begin{equation}
      \label{eq:SolODE}
      \dot p
      =
      \phi \left( t, p, A\left( \rho(t) \right) (p) \right) \,;
    \end{equation}
  \item $\rho(0) = \bar \rho$ and $p(0) = \bar p$.
  \end{enumerate}
\end{definition}
For the standard definition of Kru\v zkov solution we refer
to~\cite[Definition~1]{Kruzkov}, for that of Caratheodory solution,
see~\cite[\S~1]{Filippov}.

\begin{theorem}
  \label{thm:main}
  Under conditions~\textbf{(f)}, \textbf{($\boldsymbol{\phi}$)}
  and~\textbf{(A)}, for any initial datum $\bar p \in \reali^{N_p}$
  and $\bar \rho \in (\L1 \cap \BV) ( \reali^{N_x}; [0, R])$,
  problem~(\ref{eq:Problem}) admits a unique solution in the sense of
  Definition~\ref{def:sol}. This solution can be extended to all $I$.

  Let now $f_1$, $f_2$ satisfy~\textbf{(f)}; $A_1$, $A_2$
  satisfy~\textbf{(A)} and $\phi_1$, $\phi_2$
  satisfy~\textbf{($\boldsymbol{\phi}$)}; in all cases for the same
  interval $I$ and the same parameters or functions $R, C_f, C_A,
  C_\phi$. Then, for any initial data $(\bar\rho_1, \bar p_1),
  (\bar\rho_2, \bar p_2) \in (\L1 \cap \BV) ( \reali^{N_x}; [0, R])
  \times \reali^{N_p}$, the solutions $(\rho_1, p_1)$ and $(\rho_2,
  p_2)$ to the problems
  \begin{equation}
    \label{eq:Problem2}
    \left\{
      \begin{array}{l@{}}
        \partial_t \rho_1
        +
        \div_x f_1 \left( t, x, \rho_1, p_1(t) \right)
        =
        0
        \\
        \dot p_1
        =
        \phi_1\left( t, p_1, \left(A_1\rho_1(t) \right) (p_1) \right)
        \\
        \rho_1(0,x) = \bar \rho_1(x)
        \\
        p_1(0) = \bar p_1
      \end{array}
    \right.
    \quad \mbox{ and } \quad
    \left\{
      \begin{array}{l}
        \partial_t \rho_2
        +
        \div_x f_2 \left( t, x, \rho_2, p_2(t) \right)
        =
        0
        \\
        \dot p_2
        =
        \phi_2\left( t, p_2, \left(A_2\rho_2(t) \right) (p_2) \right)
        \\
        \rho_2(0,x) = \bar \rho_2(x)
        \\
        p_2(0) = \bar p_2
      \end{array}
    \right.
  \end{equation}
  satisfy the inequalities
  \begin{eqnarray*}
    & &
    \norma{(\rho_1-\rho_2)(t)}_{\L1}
    \\
    & \leq &
    \left( 1 + \mathcal{K}(t) \right)
    \norma{\bar \rho_1 - \bar \rho_2}_{\L1}
    \\
    & &
    +
    \mathcal{K}(t)
    \left(
      \norma{\partial_\rho(f_1-f_2)}_{\L\infty(\Omega_t)}
      +
      \norma{\div(f_1-f_2)}_{\L1(\reali^{N_x})\times \L\infty([0,t] \times
        [0,R]\times K_t)}
    \right)
    \\
    & &
    +
    \mathcal{K}(t)
    \left(
      \norma{\phi_1-\phi_2}_{\L\infty([0,t] \times K_t \times [0,C_A])}
      +
      \norma{A_1 - A_2}_{\mathcal{L}(\L1,\W{1}{\infty})}
      +
      \norma{\bar p_1 - \bar p_2}
    \right)
    \\
    \mbox{ and}
    & &
    \\
    & &
    \norma{(p_1-p_2)(t)}
    \\
    & \leq &
    \left( 1 + \mathcal{K}(t) \right)
    \norma{\bar p_1 - \bar p_2}
    \\
    & &
    +
    \mathcal{K}(t)
    \left(
      \norma{\partial_\rho(f_1-f_2)}_{\L\infty(\Omega_t)}
      +
      \norma{\div(f_1-f_2)}_{\L1(\reali^{N_x})\times \L\infty([0,t] \times
        [0,R]\times K_t)}
    \right)
    \\
    & &
    +
    \mathcal{K}(t)
    \left(
      \norma{\phi_1-\phi_2}_{\L\infty([0,t] \times K_t \times [0,C_A])}
      +
      \norma{A_1 - A_2}_{\mathcal{L}(\L1,\W{1}{\infty})}
      +
      \norma{\bar \rho_1 - \bar \rho_2}_{\L1}
    \right)
  \end{eqnarray*}
  where $\mathcal{K} \in \C0(I;\reali^+)$ vanishes at $t=0$.
\end{theorem}

\noindent More detailed expressions of the various coefficients are
presented in Section~\ref{sec:Tech}.  The proof, which is deferred to
Section~\ref{sec:Tech}, is obtained through Banach Contraction
Theorem. The necessary estimates for the convergence are a consequence
of~\cite[Theorem~5]{Kruzkov}, \cite[Theorem~2.5]{ColomboMercierRosini}
and of an adaptation of the standard theory of Caratheodory ordinary
differential equations, collected in the following two lemmas.

\begin{lemma}
  \label{lem:pde}
  Let~\textbf{(f)} hold. Choose any $\bar \rho \in (\L1 \cap \L\infty
  \cap \BV) (\reali^{N_x};[0,R])$.  Fix a function $\pi \in \C0 (I;
  \reali^{N_p})$. Then, the conservation law
  \begin{equation}
    \label{eq:HCL}
    \left\{
      \begin{array}{l}
        \partial_t \rho
        +
        \div_x f\left(t, x,  \rho, \pi(t) \right)
        =
        0
        \\
        \rho(0,x) = \bar \rho(x)
      \end{array}
    \right.
  \end{equation}
  admits a unique solution $\rho \in \C0 \left(I; \L1 (\reali^{N_x},
    [0,R])\right)$. For all $t \in I$, introduce the compact set $K_t
  = B_{\reali^{N_p}}(0, \norma{\pi}_{\C0([0,t])})$, denote $\Omega_t =
  [0,t] \times \reali^{N_x} \times [0,R] \times K_t$ and define
  \begin{equation}
    \label{eq:kappa}
    \kappa_t
    =
    (2N_x+1)
    \norma{\nabla_x\partial_\rho f}_{\L\infty(\Omega_t)}\,.
  \end{equation}
  Then, the following $\BV$ estimate holds: for all $t \in I$
  \begin{equation}
    \label{eq:TV}
    \tv \left( \rho(t) \right)
    \leq
    \left(
      \tv(\bar\rho)
      +
      N_x W_{N_x}\, t \int_{\reali^{N_x}}\norma{\nabla_x \div_x f(\cdot, x,
        \cdot, \cdot)}_{\L\infty([0,t] \times [0,R]\times K_t)} \d{x}
    \right)
    e^{\kappa_t t} \,.
  \end{equation}
  Let now, for $i = 1, 2$, $\rho_i$ be the solution
  to~(\ref{eq:Problem}) corresponding to the initial datum
  $\bar\rho_i$ and to the equation defined by $\pi_i \in
  \C0(I;\reali^{N_p})$ and by $f_i$, satisfying~\textbf{(f)}. Then,
  \begin{equation}
    \label{eq:dep_qup}
    \begin{array}{rcl}
      \norma{(\rho_1-\rho_2)(t)}_{\L1}
      & \leq &
      \norma{\bar\rho_1 - \bar\rho_2}_{\L1}
      \\
      & & +
      t \, \mathcal{C}(t)
      \Bigl[
      \norma{\pi_1-\pi_2}_{\L\infty([0,t])} +
      \norma{\partial_\rho(f_1-f_2)}_{\L\infty(\Omega_t)}
      \\
      & &
      \qquad\qquad
      +
      \norma{\div(f_1 - f_2)}_{\L1(\reali^{N_x}) \times \L\infty([0,t]
        \times [0,R] \times K_t)}
      \Bigr]
    \end{array}
  \end{equation}
  where $\mathcal{C}(t)$ depends on $\tv(\bar\rho_1)$,
  $\norma{\nabla_x \partial_\rho f_1}_{\L\infty(\Omega_t)}$,
  $\norma{\nabla_x\div_x
    f_1}_{\L1(\reali^{N_x})\times\L\infty([0,t]\times [0,R]\times
    K_t)}$ and $\norma{\nabla_p \partial_\rho
    f_2}_{\L\infty(\Omega_t)}$, $\norma{\div_x \nabla_p
    f_2}_{\L1(\reali^{N_x})\times\L\infty([0,t]\times[0,R]\times
    K_t)}$, $t$.
\end{lemma}

\noindent An explicit expression of $\mathcal{C}(t)$ is provided
in~(\ref{eq:explicit}).

The estimates related to the ordinary differential equation are
provided by the following lemma.

\begin{lemma}
  \label{lem:edo}
  Let~\textbf{($\boldsymbol{\phi}$)} and~\textbf{($\boldsymbol{A}$)}
  hold. Choose an initial datum $\bar p \in \reali^{N_p}$ and fix a
  function $r\in \C0\left(I; \L1(\reali^{N_x};[0,R]) \right)$. Then,
  the ordinary differential equation
  \begin{equation}
    \label{eq:edo}
    \left\{
      \begin{array}{l}
        \dot p
        =
        \phi\left( t, p, A\left(r(t)\right)(p) \right)
        \\
        p(0) =  \bar p\,.
      \end{array}
    \right.
  \end{equation}
  admits a unique solution $p\in \Wloc{1}{\infty}(I;
  \reali^{N_p})$. The following bound holds:
  \begin{equation}
    \label{eq:SupBound}
    \norma{p(t)}
    \leq
    \left(\norma{\bar p} + 1 \right)
    e^{\int_0^t C_\phi(\tau) \d\tau}
    - 1
    \,.
  \end{equation}
  Given two initial conditions $\bar p_1, \bar p_2 \in \reali^{N_p}$,
  two functions $r_1, r_2\in \C0\left(I; \L1(\reali^{N_x};[0,R])
  \right)$, two speed laws $\phi_1, \phi_2$
  satisfying~\textbf{($\boldsymbol{\phi}$)} and two averaging
  operators $A_1$, $A_2$ satisfying~\textbf{(A)}, define
  \begin{equation}
    \label{eq:F}
    F(t)
    =
    \left(
      1
      +
      C_A
      \norma{r_1}_{\L\infty([0,t];\L1)}
    \right)
    \int_0^t C_\phi(\tau) \d\tau \,.
  \end{equation}
  Then,
  \begin{equation}
    \label{eq:est_edo}
    \begin{array}{@{}r@{\;}c@{\;}l@{}}
      & &
      \displaystyle
      \norma{(p_1-p_2)(t)}
      \\
      & \leq &
      \displaystyle
      e^{F(t)} \norma{\bar p_1 - \bar p_2}
      +
      \int_0^t e^{F(t)-F(\tau)}
      \norma{\phi_1(\tau) - \phi_2(\tau)}_{\L\infty} \d\tau
      \\
      & &
      +
      \displaystyle
      \int_0^t \! e^{F(t)-F(\tau)} C_\phi(\tau)
      \left(
        C_A
        \norma{(r_1-r_2)(\tau)}_{\L1}
        +
        \norma{A_1-A_2}_{\mathcal{L}(\L1,\W{1}{\infty})}
        \norma{r_2(\tau)}_{\L1}
      \right)
      \d\tau .\!\!\!\!
    \end{array}
  \end{equation}
\end{lemma}

In the applications below, the support of the initial data is
compact. Thanks to the finite propagation speed typical of
conservation laws, this allows a major simplification in the
assumptions of Theorem~\ref{thm:main}.

\begin{corollary}
  \label{cor:compact}
  Consider problem~\eqref{eq:Problem} with $f$
  satisfying~\textbf{(f.\ref{it:f1})}, \textbf{(f.\ref{it:f2})}
  and~\textbf{(f.\ref{it:f3})}. Let $A$ satisfy~\textbf{(A)} and
  $\phi$ satisfy~\textbf{($\boldsymbol{\phi}$)}. If $\bar\rho$
  vanishes outside a compact set, then problem~\eqref{eq:Problem}
  admits a unique solution in the sense of
  Definition~\ref{def:sol}. This solution can be extended to all of
  $I$. Moreover, the stability estimates of Theorem~\ref{thm:main}
  apply, provided both $\bar\rho_1$ and $\bar \rho_2$ vanish outside a
  compact set.
\end{corollary}

\section{Applications}
\label{sec:Appl}

This section is devoted to a few sample applications of
Theorem~\ref{thm:main}. While the unknown $\rho$ keeps throughout the
meaning of a scalar density, the state $p$ of the individuals is the
position of a single agent in~\S~\ref{subs:PP}, it is a vector of
several positions in~\S~\ref{subs:Dog} and it becomes a 4--vector
position--speed in~\S~\ref{subs:Hawk}.

Numerical integrations are also provided in order to show the
qualitative behavior of the solutions. In all cases, the
Lax--Friedrichs method, see~\cite[\S~12.5]{LeVequeBook}, with
dimensional splitting was used for the conservation law and Euler
polygonals to integrate the ordinary differential equation.

\subsection{The Pied Piper}
\label{subs:PP}

\begin{figure}[htpb]
  \centering
  \begin{psfrags}
    \psfrag{0}{$0$} \psfrag{v}{$v$} \psfrag{fun}{$\rho \to v(\rho)$}
    \psfrag{R}{$R$} \psfrag{V}{ } \psfrag{rho}{$\rho$}
    \includegraphics[width=45mm]{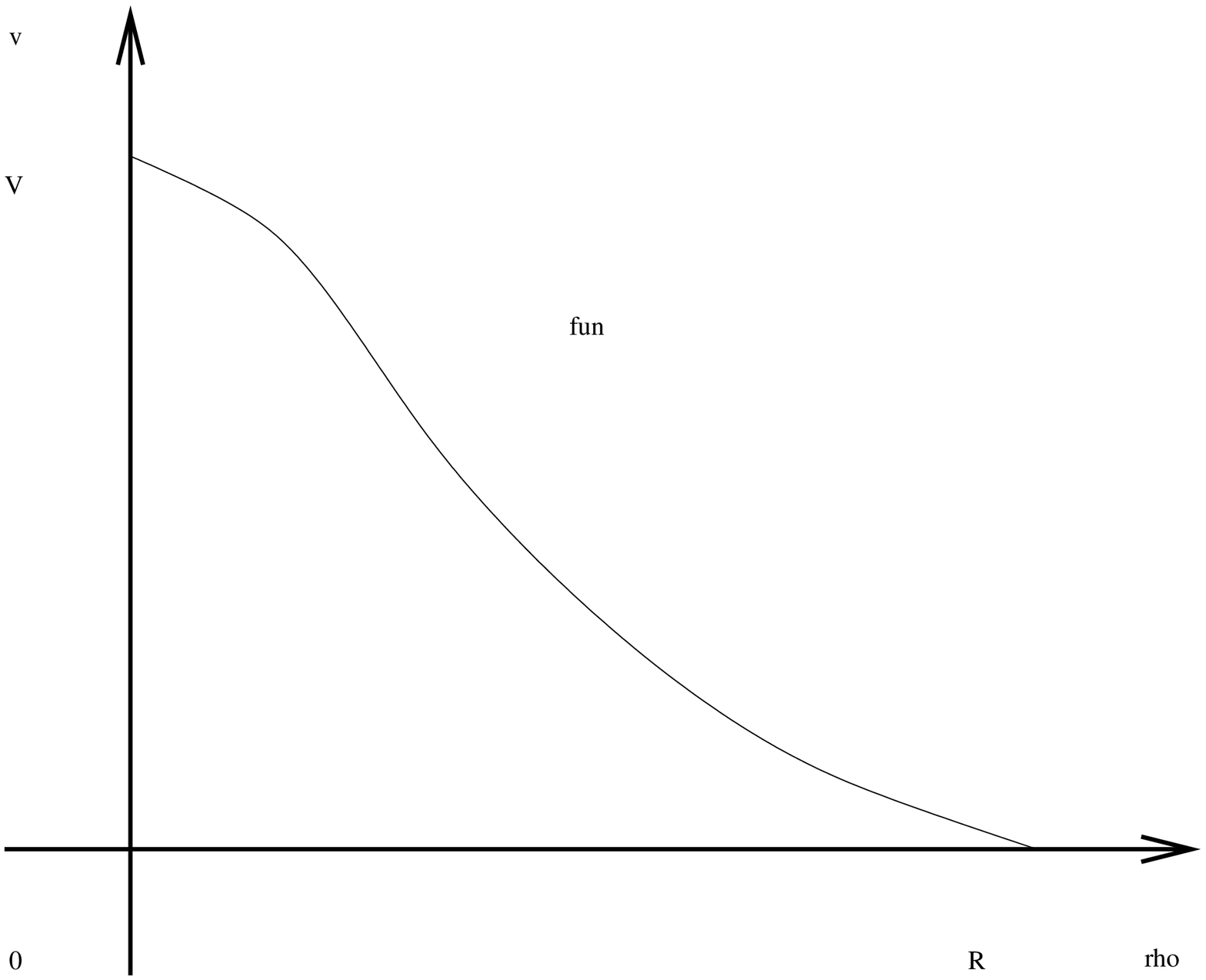}
  \end{psfrags}
  \hfil
  \begin{psfrags}
    \psfrag{0}{$0$} \psfrag{x}{$\norma{x}$}
    \psfrag{mv}{$\norma{\vec{\boldsymbol{v}}}$} \psfrag{fun}{$x \to
      \vec{\boldsymbol{v}}(x)$}
    \includegraphics[width=45mm]{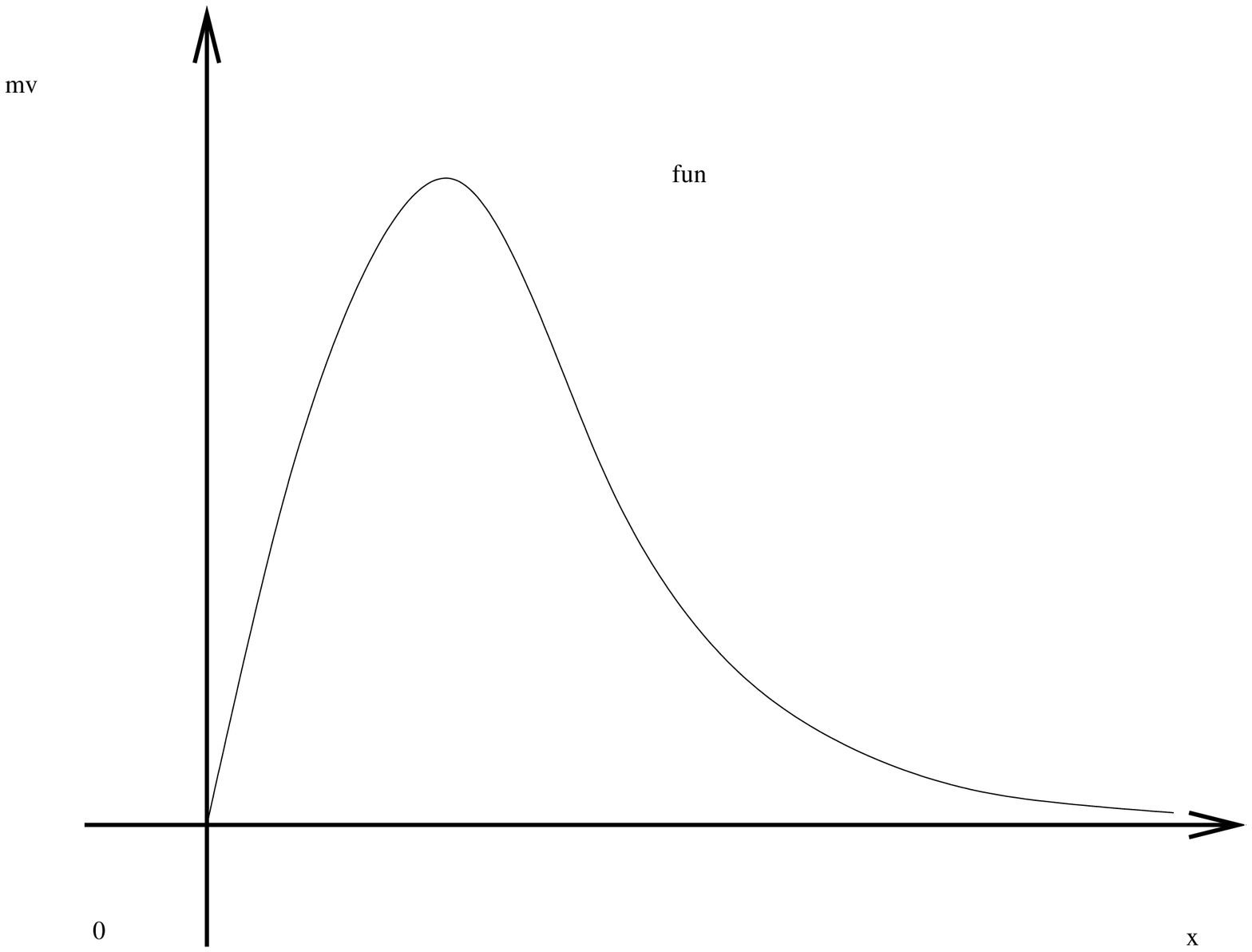}
  \end{psfrags}
  \hfil
  \begin{psfrags}
    \psfrag{0}{$0$} \psfrag{v}{$q$} \psfrag{fun}{$\rho \to q(\rho)$}
    \psfrag{R}{$R$} \psfrag{V}{ } \psfrag{rho}{$\rho$}
    \includegraphics[width=45mm]{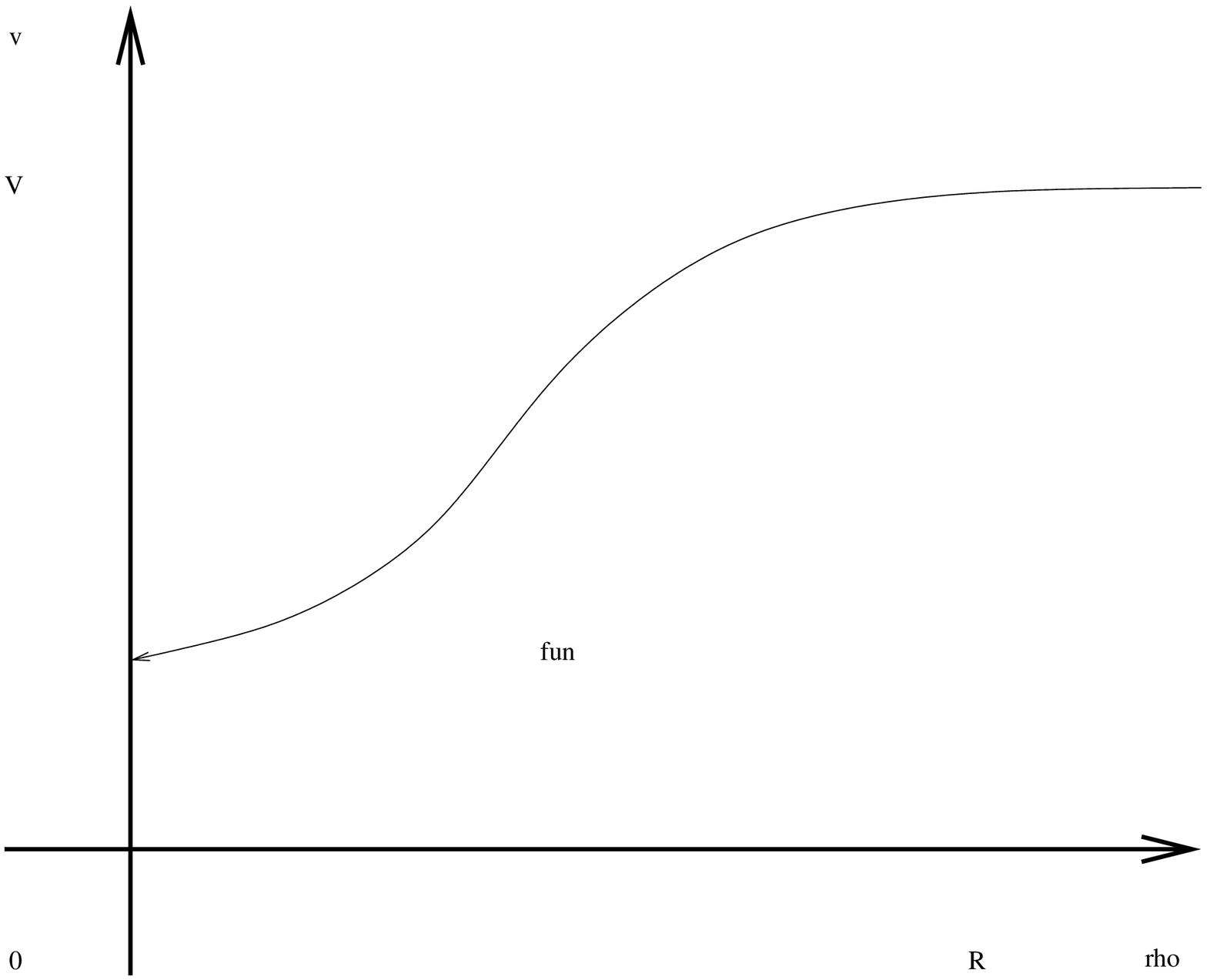}
  \end{psfrags}
  \Caption{Left, $v$ is assumed $\C2$ and decreasing. Center,
    $\vec{\boldsymbol{v}}$ describes the attraction felt by the mice
    towards the piper. Right, $q$ accounts for the acceleration of the
    piper when surrounded by a high mice density.}
  \label{fig:ppAssumptions}
\end{figure}

As a first toy application we consider the situation described
in~\cite[n.~246]{GrimmBrothers}. To lure rats away, the city of
Hamelin (now Hamel) hires a rat-catcher who, playing his magic pipe,
attracts all mice out of the city. In this case, $\rho = \rho(t,x)$ is
the mice density and $p = p(t)$ is the position of the piper. Rats
move with a speed $v(\rho) \, \vec{\boldsymbol{v}} (p-x)$, with the
scalar $v$ and the vector $\vec{\boldsymbol{v}}$ having the
qualitative behavior in Figure~\ref{fig:ppAssumptions}.
More precisely, at density $0$ mice have the fastest speed while at
density $R$ their speed vanishes. The term $\vec{\boldsymbol{v}}$
accounts for the attraction of the mice towards the piper. The magic
musician has a speed $q(\rho*\eta)\, \vec{\boldsymbol{\psi}}(t)$,
i.e.~he moves faster when the average density of mice around him is
higher. On the contrary, when only few rats are near to him, he slows
down.

\begin{lemma}
  Let $N_x = 2$, $N_p = 2$, $N_r = 1$ and fix a positive $R$. Assume
  $v \in \C2([0,R]; \reali)$, $\vec{\boldsymbol{v}} \in
  \C2(\reali^2;\reali^2)$, $q \in \W1\infty([0,R];\reali)$,
  $\vec{\boldsymbol{\psi}} \in \W1\infty(\reali^+;\reali^2)$, $\eta
  \in \Cc2(\reali^2, \reali)$ with $\int_{\reali^2} \eta \d{x}
  =1$. Assume that $v(R) = 0$.  Define
  \begin{equation}
    \label{eq:pp_f}
    f(t,x,\rho,p)
    =
    \rho \, v(\rho) \, \vec{\boldsymbol{v}} (p-x)
    \qquad
    \phi(t,p,r)
    =
    q(r) \, \vec{\boldsymbol{\psi}}(t)
    \qquad
    A\rho
    =
    \rho *_x \eta \,.
  \end{equation}
  Then, this setting fits in the framework of
  Corollary~\ref{cor:compact} as soon as $\bar \rho$ vanishes outside
  a compact set.
\end{lemma}

The proof is immediate and, hence, omitted.

\smallskip

\textbf{Numerical example:} To fix a specific situation, we choose the
following functions in~(\ref{eq:Problem}):
\begin{equation}
  \label{eq:PiedPiper}
  \begin{array}{rcl@{\qquad}l@{\qquad}l}
    v(\rho)
    & = &
    V_{\max}\, \left(1-\frac{\rho}{R}\right)
    &
    V_{\max} = 9
    &
    R = 1
    \\
    \vec{\boldsymbol{v}} (x)
    & = &
    x \, e^{-\norma{x}^2}
    \\
    q(r)
    & = &
    v_p + \frac{V_p-v_p}{R} r
    &
    V_p = 7
    &
    v_p = 1
    \\
    \vec{\boldsymbol{\psi}}(t)
    & = &
    \left[
      \begin{array}{c}
        \cos \omega t \\ - \sin \omega t
      \end{array}
    \right]
    &
    \omega = 1
    \\
    \eta(x)
    & = &
    \frac{3}{\pi {r_p}^6}\,
    \left( \max \left\{0, {r_p}^2 - \norma{x}^2 \right\} \right)^2
    &
    r_p = 0.15
  \end{array}
\end{equation}
At time $t=0$, we assume that rats are uniformly distributed with
density $R=1$ in the rectangle $[-0.5, 0] \times [0.35, 0.85]$. The
piper starts moving at the point $(-1, 0.5)$.
\begin{figure}[htpb]
  \includegraphics[width=0.33\textwidth]{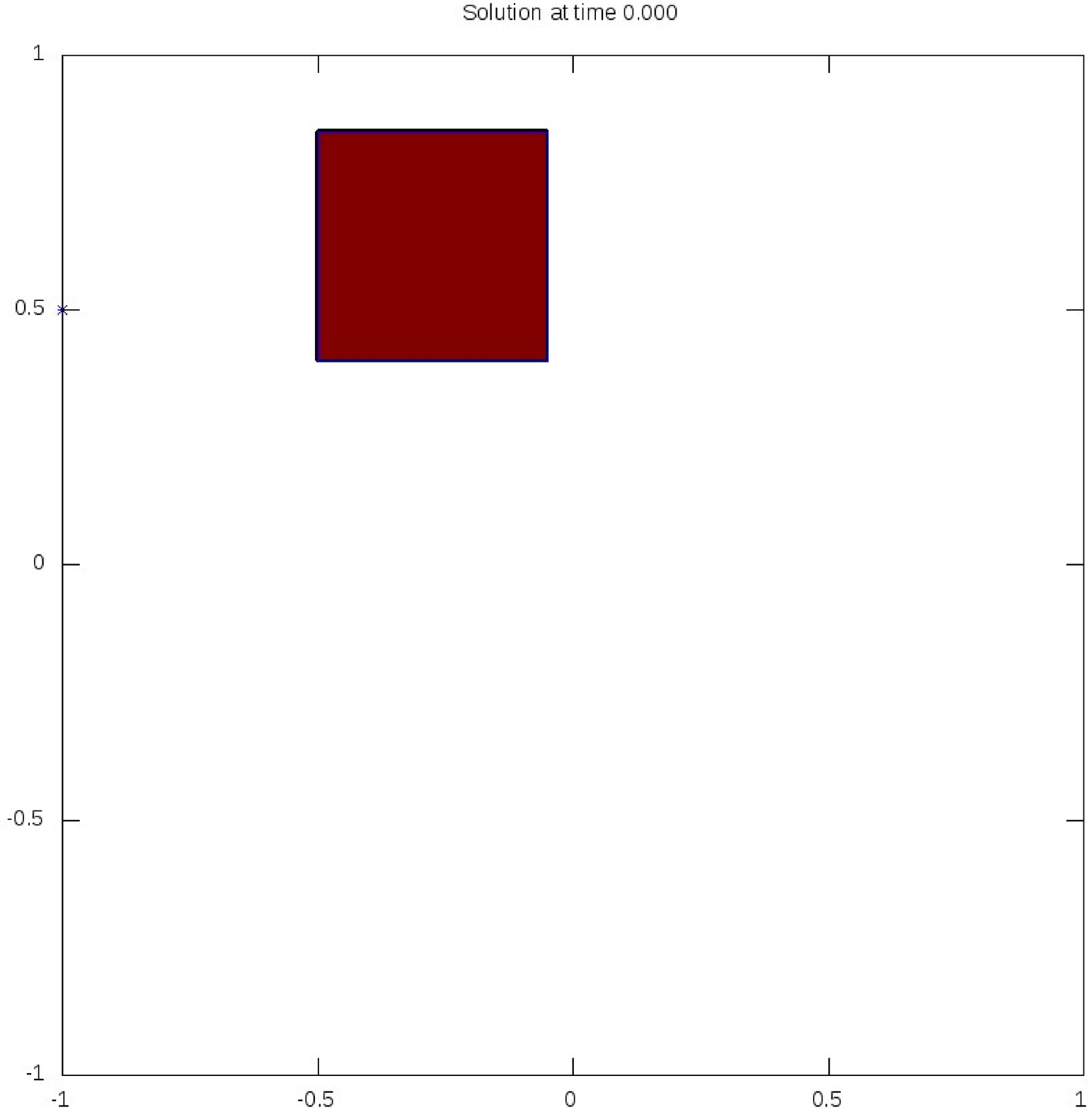}%
  \includegraphics[width=0.33\textwidth]{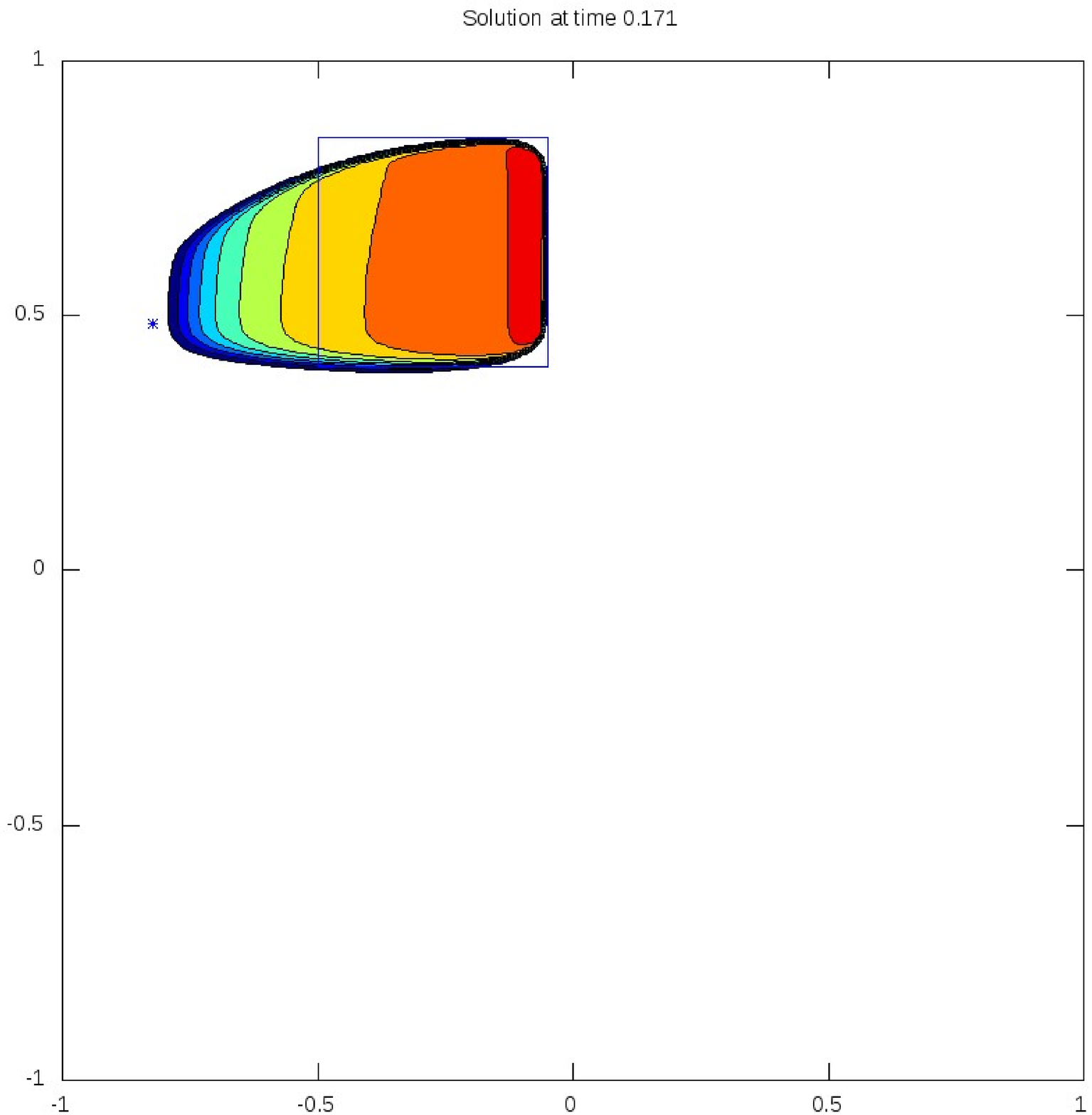}%
  \includegraphics[width=0.33\textwidth]{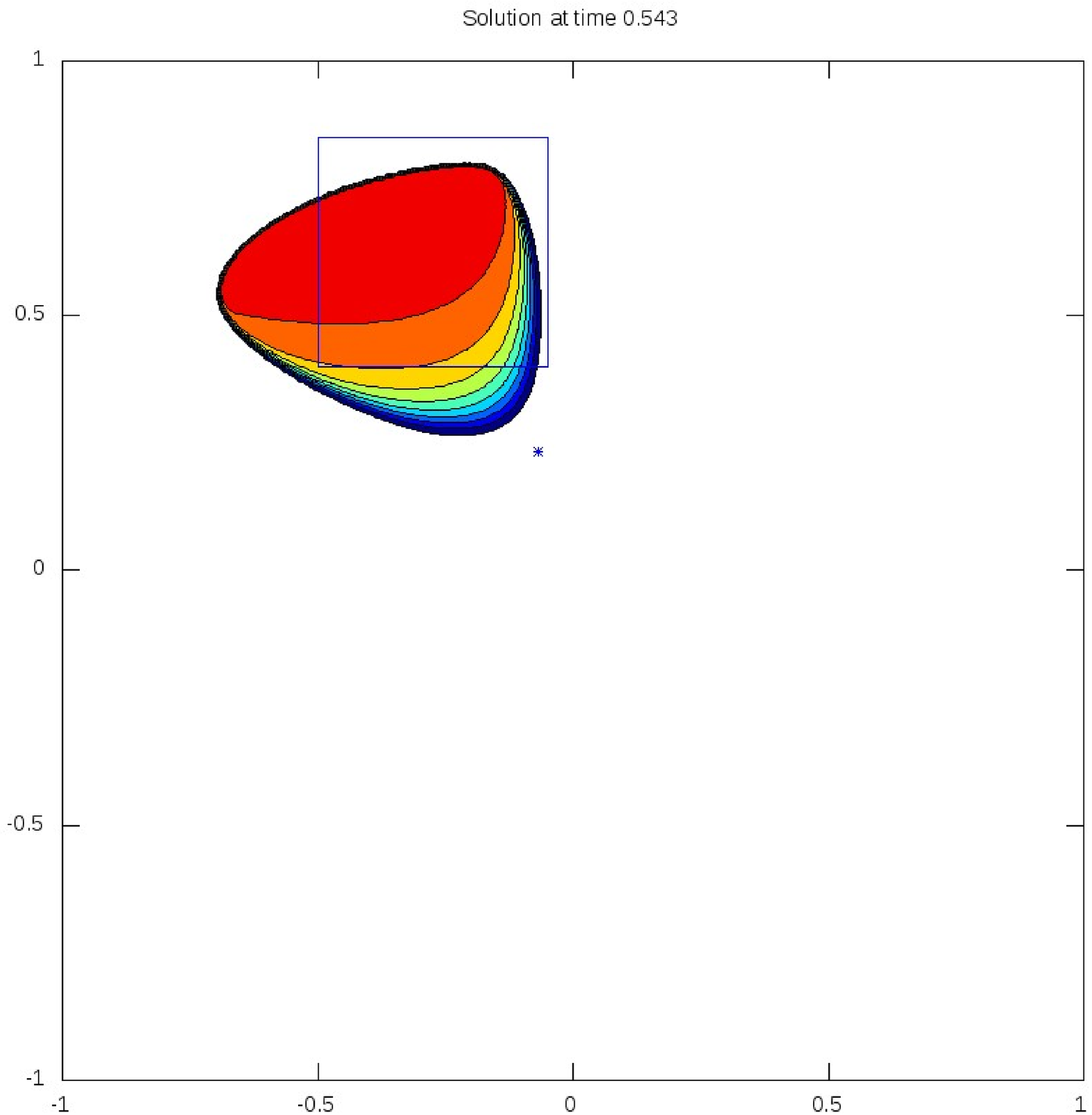}
  \\
  \includegraphics[width=0.33\textwidth]{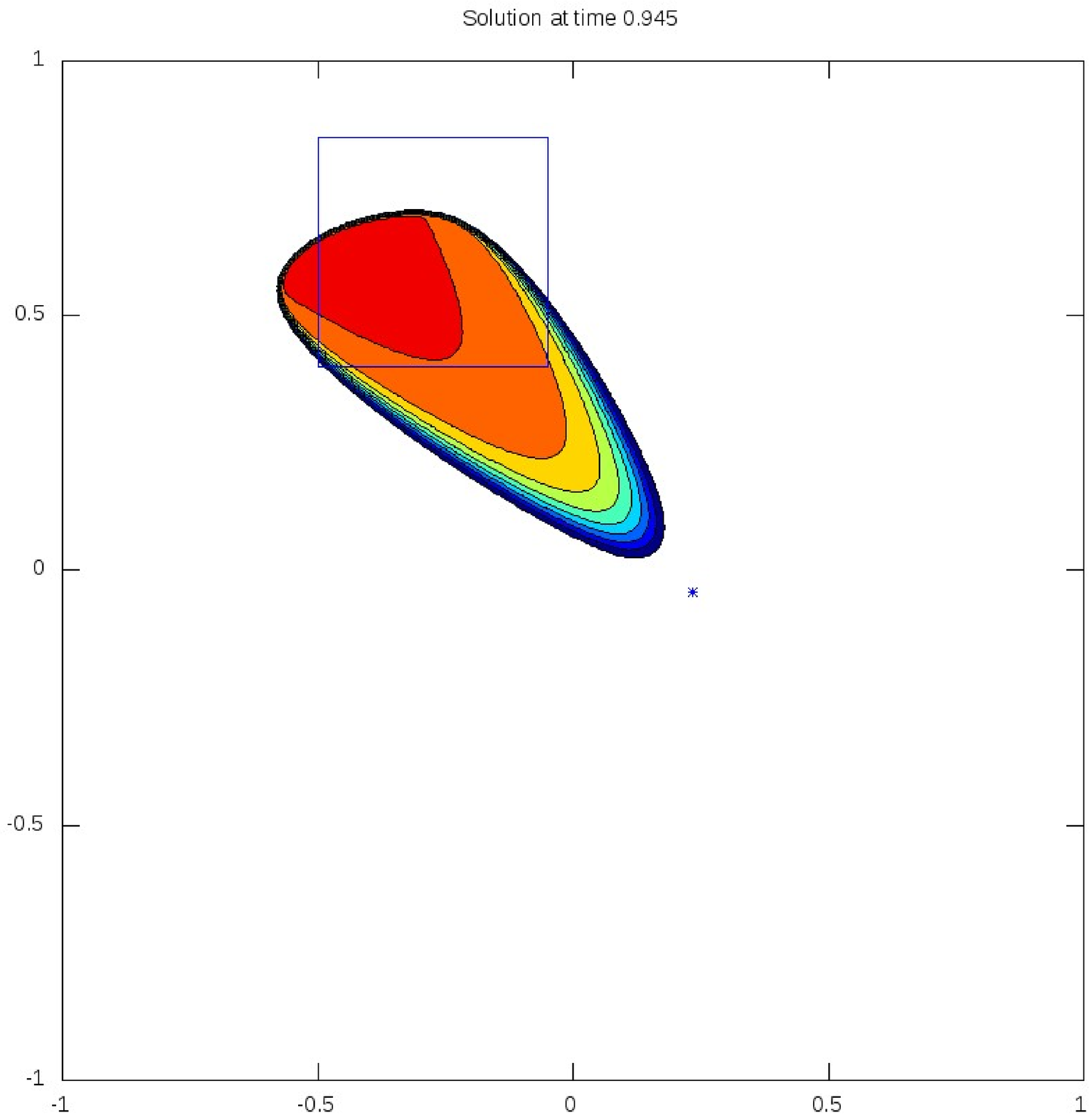}%
  \includegraphics[width=0.33\textwidth]{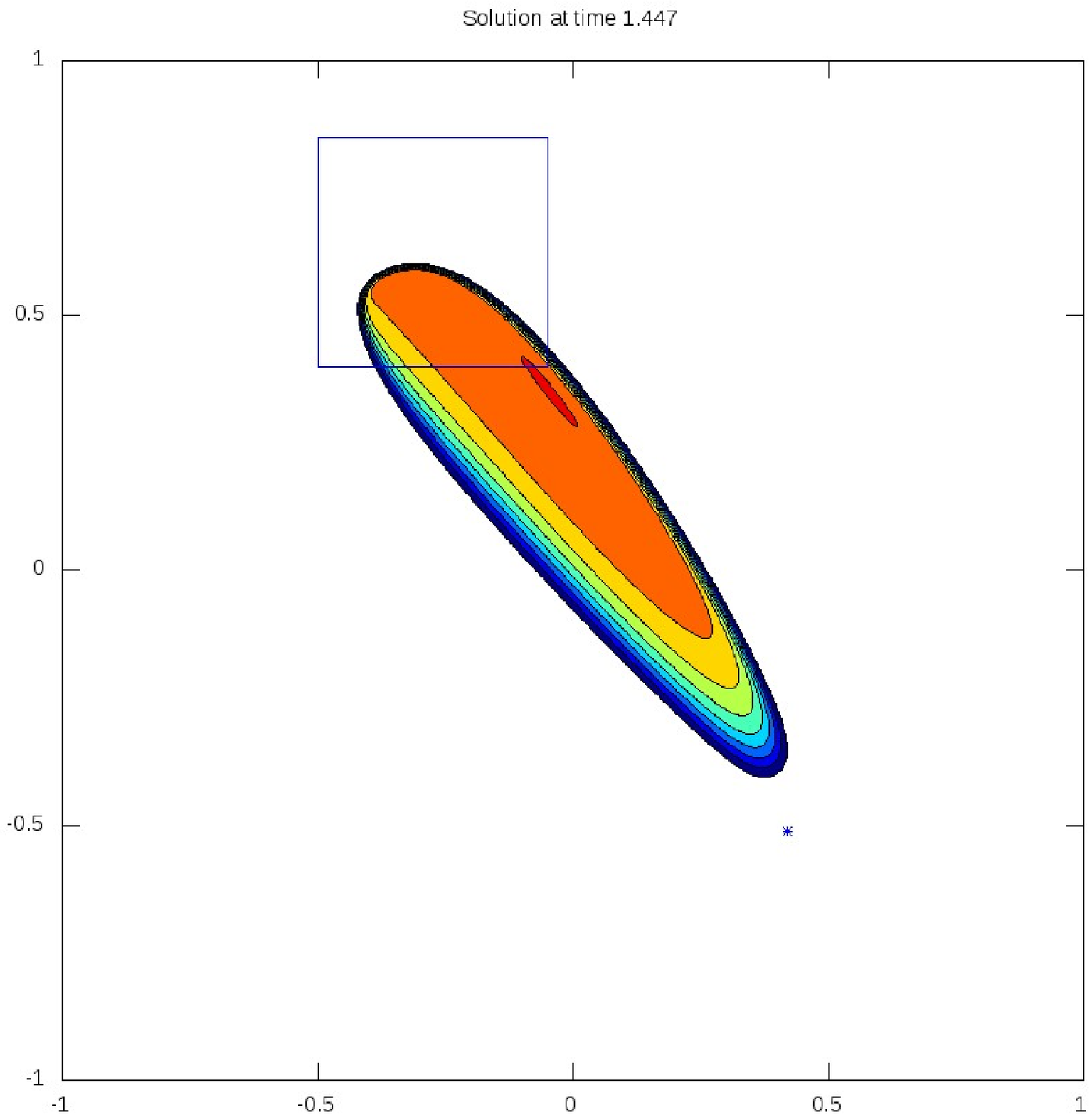}%
  \includegraphics[width=0.33\textwidth]{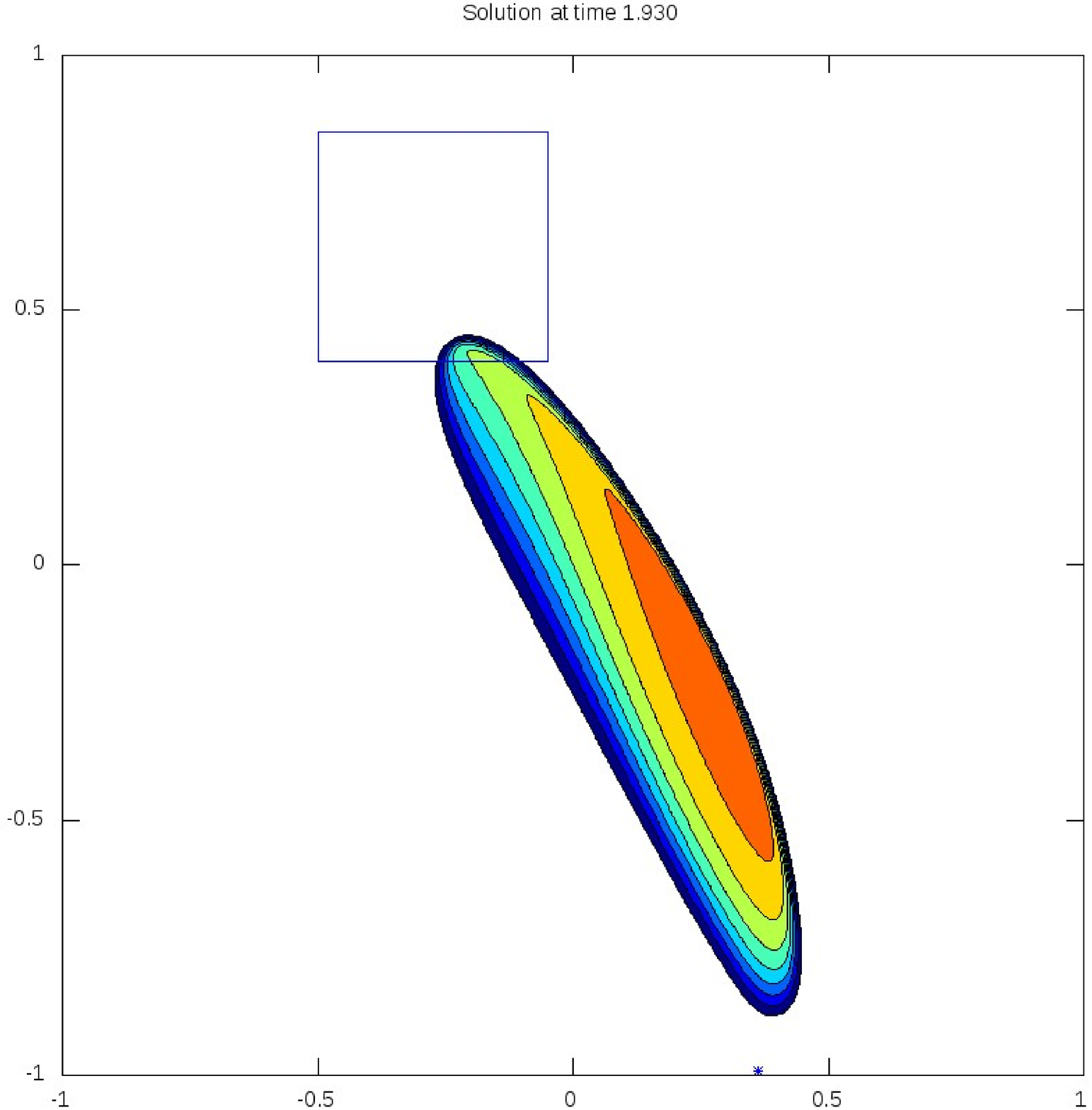}
  \Caption{The pied piper and the rats, at times $0$, when $p = (0,
    0.5)$; $0.171$, $0.543$, $0.945$, $1.447$ and $1.930$, when the
    rats almost completely left the rectangle and $p = (0.366,
    -0.983)$.}\label{fig:flute}
\end{figure}

Several optimization problems can now be stated with reference
to~(\ref{eq:Problem})--\eqref{eq:pp_f}--(\ref{eq:PiedPiper}). Referring
to the situation~\cite[n.~246]{GrimmBrothers}, a first natural
question is the following. Let the compact set $K$ be the area of the
city and fix a finite positive time $T_{\max}$. Then, find the initial
position $\bar p$ and the trajectory $\vec{\boldsymbol{\psi}}$ of the
piper so that the amount of mice left in the city at time $T_{\max}$
is minimal. In other words, we want to minimize the functional
\begin{displaymath}
  (\bar p, \vec{\boldsymbol{\psi}})
  \mapsto
  \int_K \left(\rho(\bar p, \vec{\boldsymbol{\psi}})\right)
  (T_{\max},x) \d x \,.
\end{displaymath}
Here, $\rho(\bar p, \vec{\boldsymbol{\psi}})$ is the $\rho$--component
of the solution
to~(\ref{eq:Problem})--\eqref{eq:pp_f}--(\ref{eq:PiedPiper}). The
existence of such an optimal strategy for the piper follows from
Theorem~\ref{thm:main} via a standard application of {Weierstra\ss}
Theorem.

\begin{proposition}
  \label{prop:Piper}
  Let $T_{\max}$ be finite. Denote by $K \subset \reali^2$ the compact
  Hamelin urban area. Define the set of the possible piper's route
  choices
  \begin{displaymath}
    \mathcal{K}
    =
    \left\{
      (\bar p, \vec{\boldsymbol{\psi}}) \in K \times \W{1}{\infty} (I;\reali^2)
      \colon
      \norma{\vec{\boldsymbol{\psi}}}_{\W1\infty} \leq 1
    \right\}
  \end{displaymath}
  and call $\mathcal{J} \colon \mathcal{K} \mapsto \reali$ the
  functional giving the total amount of mice in Hamelin at time
  $T_{\max}$, i.e.
  \begin{displaymath}
    \mathcal{J}(\bar p, \vec{\boldsymbol{\psi}})
    =
    \int_K \left(\rho(\bar p, \vec{\boldsymbol{\psi}})\right)
    (T_{\max},x) \d x \,,
  \end{displaymath}
  where $\rho(\bar p, \vec{\boldsymbol{\psi}})$ is the solution
  to~(\ref{eq:Problem})--(\ref{eq:pp_f})--(\ref{eq:PiedPiper}). Then,
  there exists an optimal trajectory $(\bar p_*,
  \vec{\boldsymbol{\psi}}_*) \in \mathcal{K}$ such that
  $\mathcal{J}(\bar p_*, \vec{\boldsymbol{\psi}}_*) =
  \min_{\mathcal{K}} \mathcal{J}(\bar p, \vec{\boldsymbol{\psi}})$.
\end{proposition}

\noindent Thanks to the stability estimates in Theorem~\ref{thm:main},
the proof of this proposition directly follows from Ascoli-Arzel\`a
Theorem that allows to prove the compactness of $\mathcal{K}$.

\subsection{Shepherd Dogs}
\label{subs:Dog}

On the plane, consider a herd of, say, sheeps controlled by $n$
shepherd dogs. Then, $\rho$ is the density of sheeps and $p \equiv
(p_1, \ldots, p_n)$ is the vector of the positions of the dogs, so
that each $p_i$ is in $\reali^2$. We assume that initially the sheeps
are distributed around, say, the origin and tend to disperse moving
radially with a speed directed by $\vec{\boldsymbol{v_r}}(x)$. The
duty of the dogs is to prevent this dispersion and they pursue this
goal moving around sheeps or, more precisely, with a speed $\phi$
orthogonal to the gradient of the sheeps' density. The sheeps modify
their speed escaping from the dogs with a repulsive speed
$\vec{\boldsymbol{v_d}}(x, p) = \sum_i \vec{\boldsymbol{v}}(x-p_i)$,
where $\vec{\boldsymbol{v}}$ behaves qualitatively as in
Figure~\ref{fig:ppAssumptions}. Finally, the speed of the sheeps is
then given by $v(\rho) \,( \vec{\boldsymbol{v_r}}(x) +\sum_{i=1}^n
\vec{\boldsymbol{v}} (x-p_i))$ where $v$ is maximal at the density
zero and vanishes at the maximal density $R$. This last fact means
that the sheeps can not move when their density is maximal.

\begin{lemma}
  \label{lem:Dog}
  Let $n \in \naturali$, $N_x = 2$, $N_p = 2n$, $N_r = 2n$ and fix a
  positive $R$. Assume $v \in \C2([0,R]; \reali)$,
  $\vec{\boldsymbol{v_r}}\in \C2(\reali^2;\reali^2)$
  $\vec{\boldsymbol{v}} \in \C2(\reali^2;\reali^2)$, $\eta \in
  \Cc2(\reali^2, \reali)$ with $\int_{\reali^2} \eta \d{x} =1$. Assume
  that $v(R) = 0$.  Define
  \begin{equation}
    \label{eq:dog_f}
    \begin{array}{rcl}
      f(t,x,\rho,p)
      & = &
      \rho \, v(\rho) \left( \boldsymbol{\vec{v}_r}(x)
        +
        \sum_{i=1}^n \vec{\boldsymbol{v}} (x-p_i) \right),
      \\
      \phi(t,p,r)
      & = &
      V_d\,\frac{r^\perp}{\sqrt{1+\norma{r}^2}},
      \\
      A\rho
      & = &
      \rho *_x\nabla \eta \,.
    \end{array}
  \end{equation}
  Then, this setting fits in the framework of
  Corollary~\ref{cor:compact} as soon as $\bar \rho$ vanishes outside
  a compact set.
\end{lemma}

\noindent Here, $r \equiv (r_1, \ldots, r_n)$ is a vector in
$(\reali^2)^n$ and we set $r^\perp \equiv (r_1^\perp, \ldots,
r_n^\perp)$, with $\left[
  \begin{array}{@{}c@{}}
    a\\b
  \end{array}
\right]^\perp = \left[
  \begin{array}{@{}c@{}}
    b\\-a
  \end{array}
\right]$.

\noindent\textbf{Numerical Integration:} To fix a specific situation,
we choose $n=2$ and the following functions in~(\ref{eq:Problem}):
\begin{equation}
  \label{eq:ShepDog}
  \begin{array}{rcl@{\qquad}l@{\qquad}l}
    v(\rho)
    & = &
    \displaystyle  V_{\max}\, \left(1-\frac{\rho}{R}\right)
    &
    V_{\max} = 1
    &
    R = 1
    \\
    \boldsymbol{\vec{v}} (x)
    & = &
    \displaystyle \frac{\alpha}{\sqrt{\ell}}  \,
    e^{-\norma{x}^2/\ell} \, x
    &
    \alpha = 20
    &
    \ell = 0.2
    \\
    \boldsymbol{\vec{v}_r}(x)
    &=&
    \displaystyle\frac{ \beta\, x}{1+\norma{x}^2}
    &
    \beta = 1
    \\
    \eta(x)
    & = &
    \displaystyle\frac{3}{\pi {r_p}^6}\,
    \left( \max \left\{0, {r_p}^2 - \norma{x}^2 \right\} \right)^2
    &
    r_p = 1
    \\
    V_d &=& 100
  \end{array}
\end{equation}
At time zero, sheeps are uniformly distributed at the maximal density
$R=1$ in the circumference centered at $(0,0)$ with radius $0.2$. Dogs
start moving from $(0.7, \, 0)$ and $(-0.7, \, 0)$
\begin{figure}[htpb]
  \includegraphics[width=0.33\textwidth]{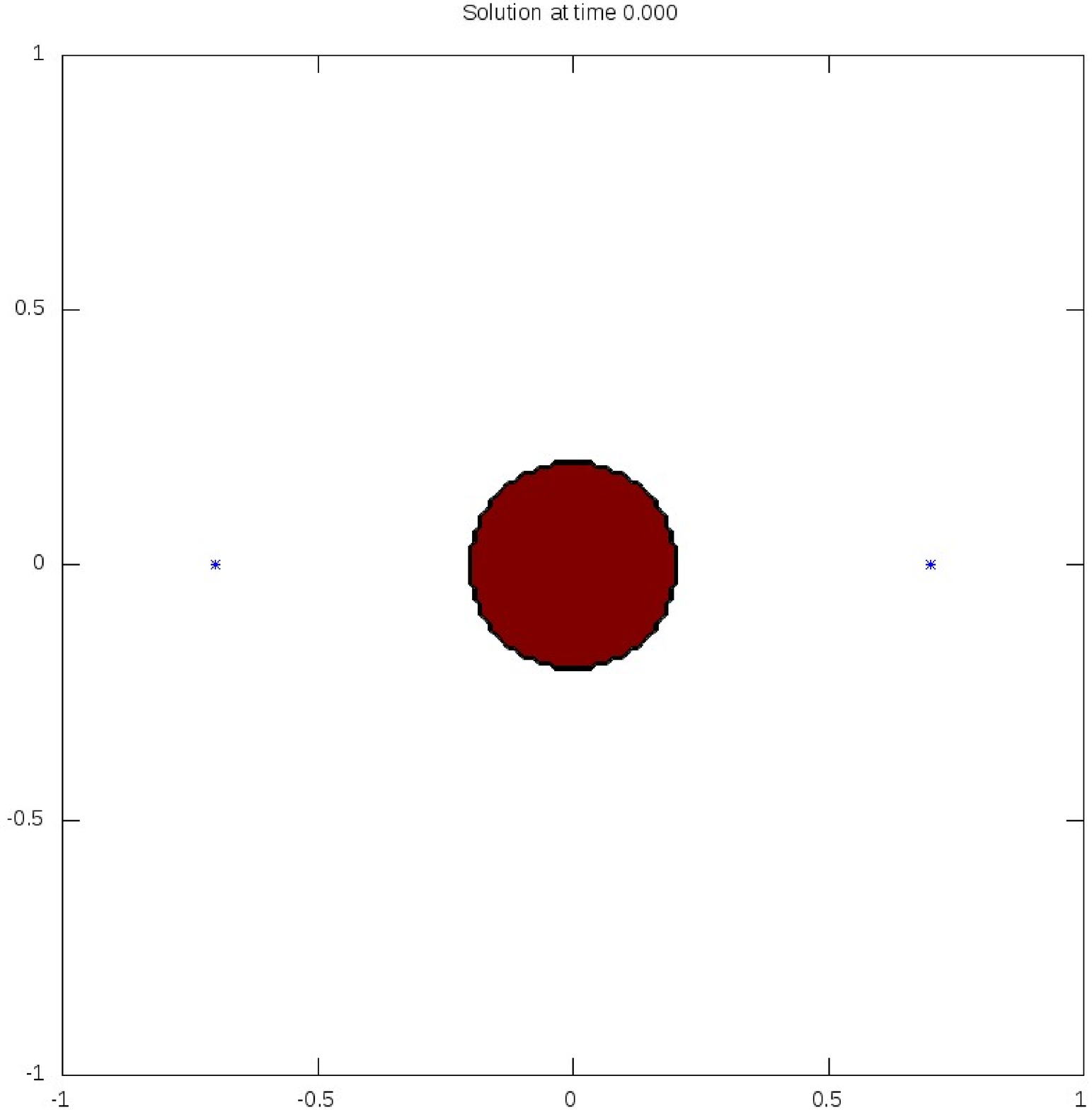}%
  \includegraphics[width=0.33\textwidth]{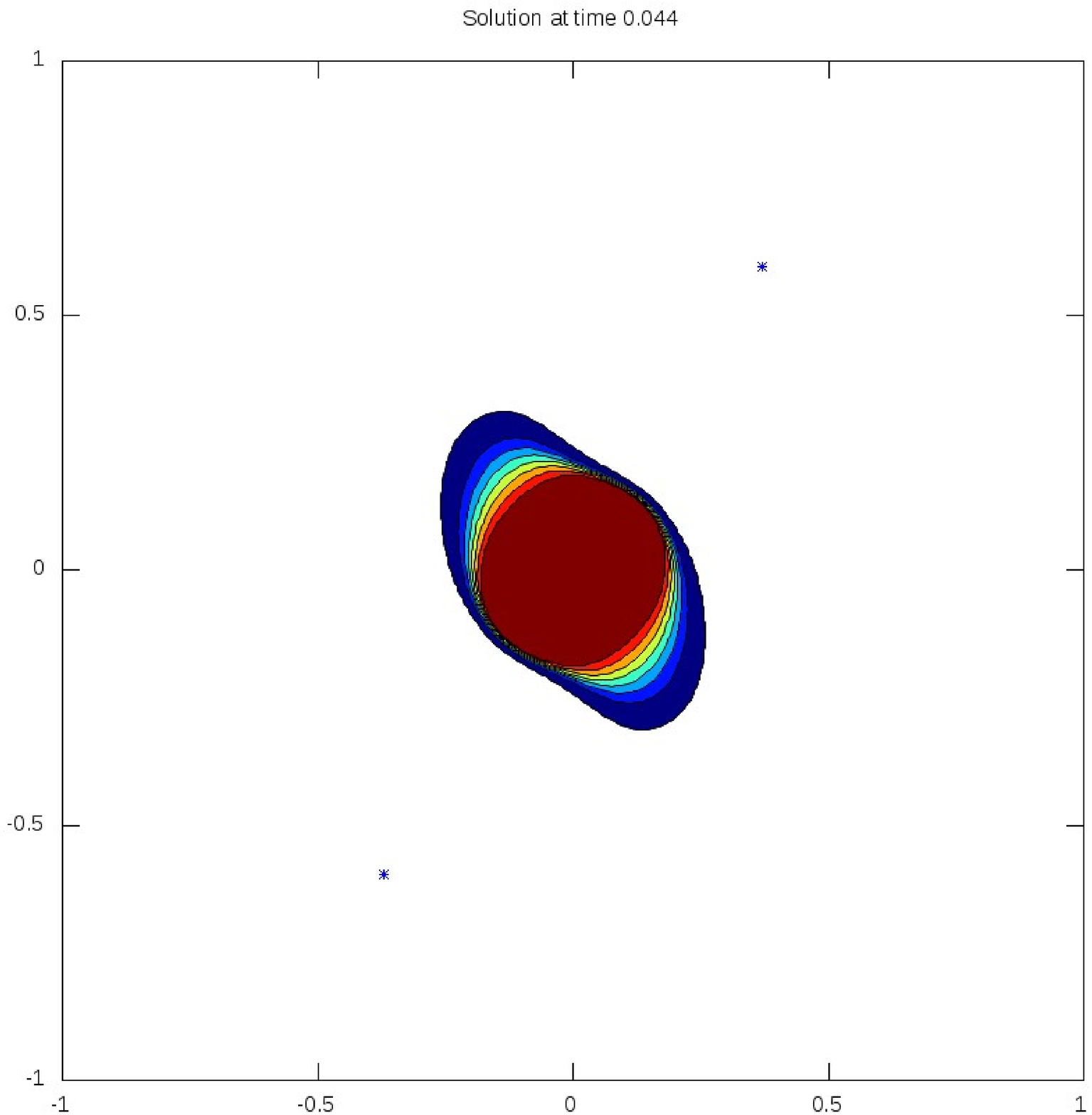}%
  \includegraphics[width=0.33\textwidth]{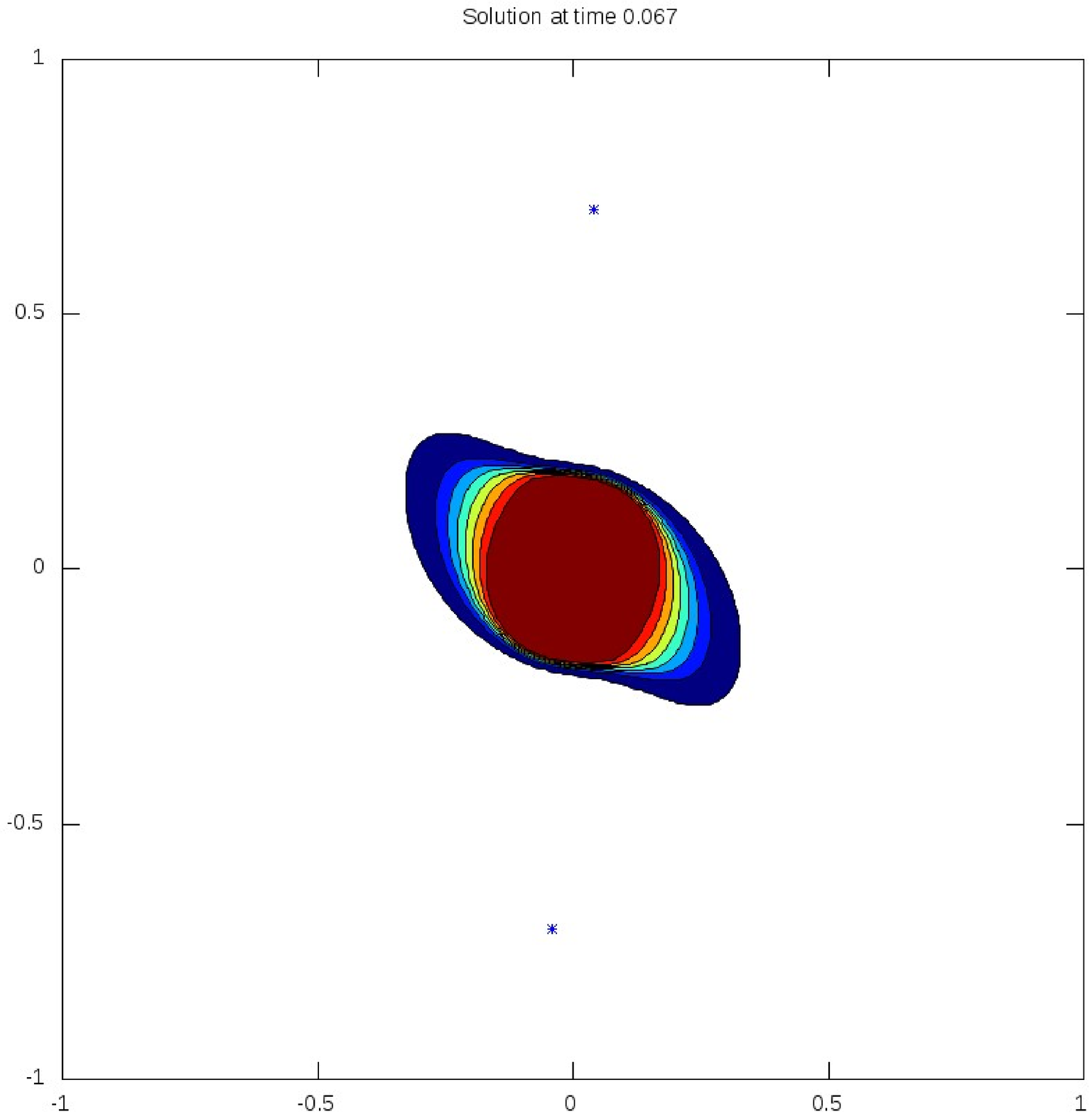}
  \\
  \includegraphics[width=0.33\textwidth]{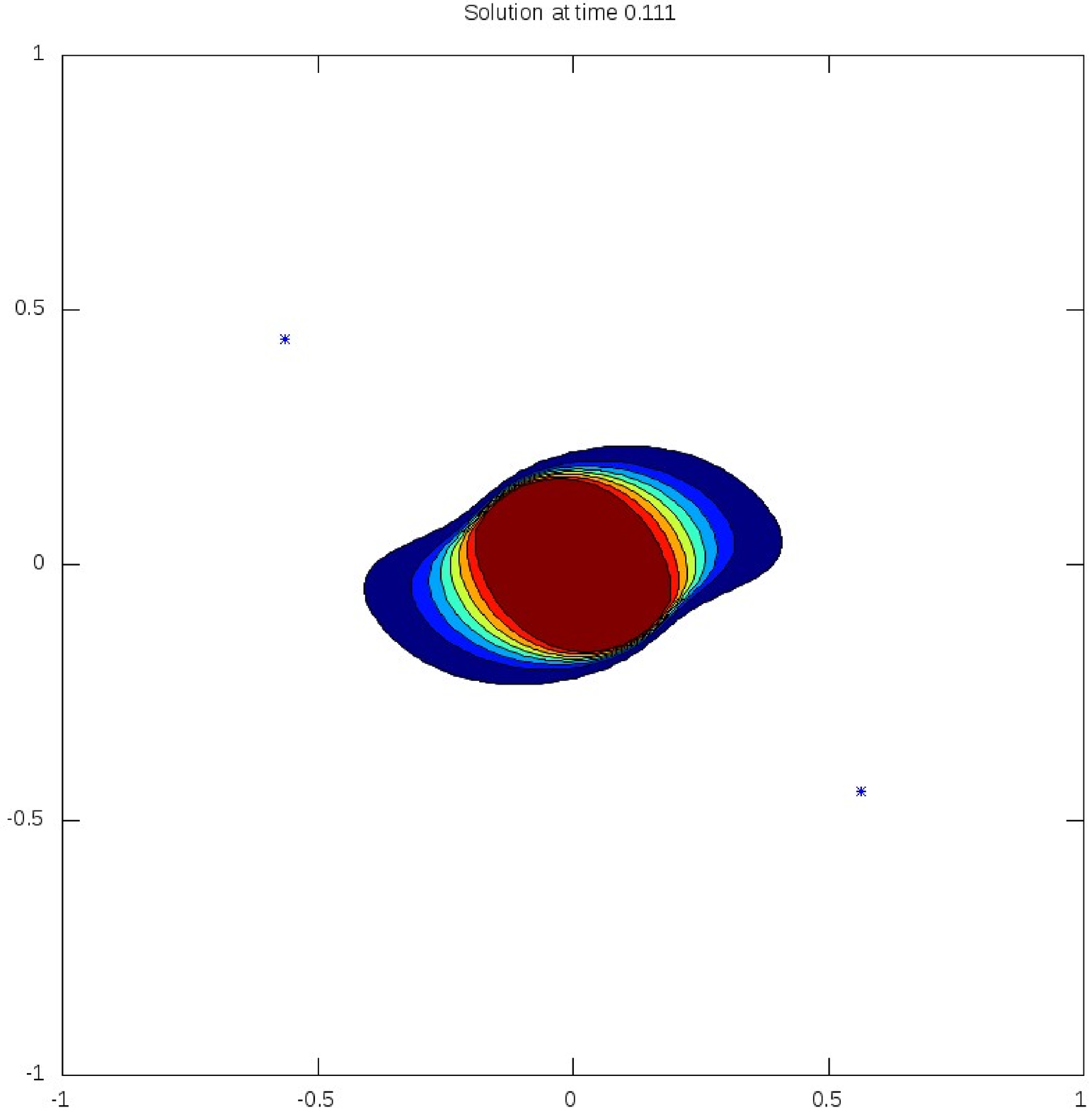}%
  \includegraphics[width=0.33\textwidth]{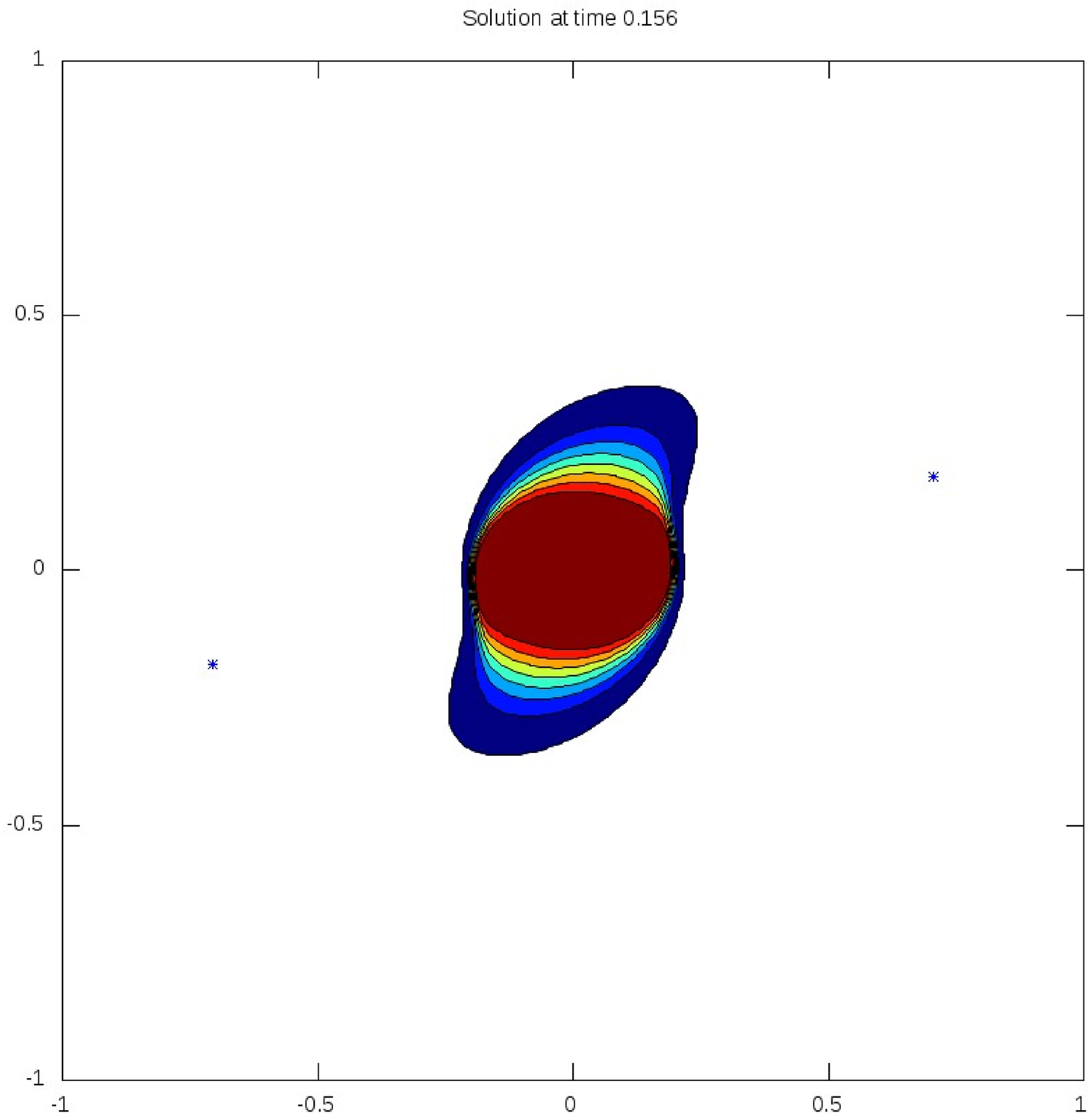}%
  \includegraphics[width=0.33\textwidth]{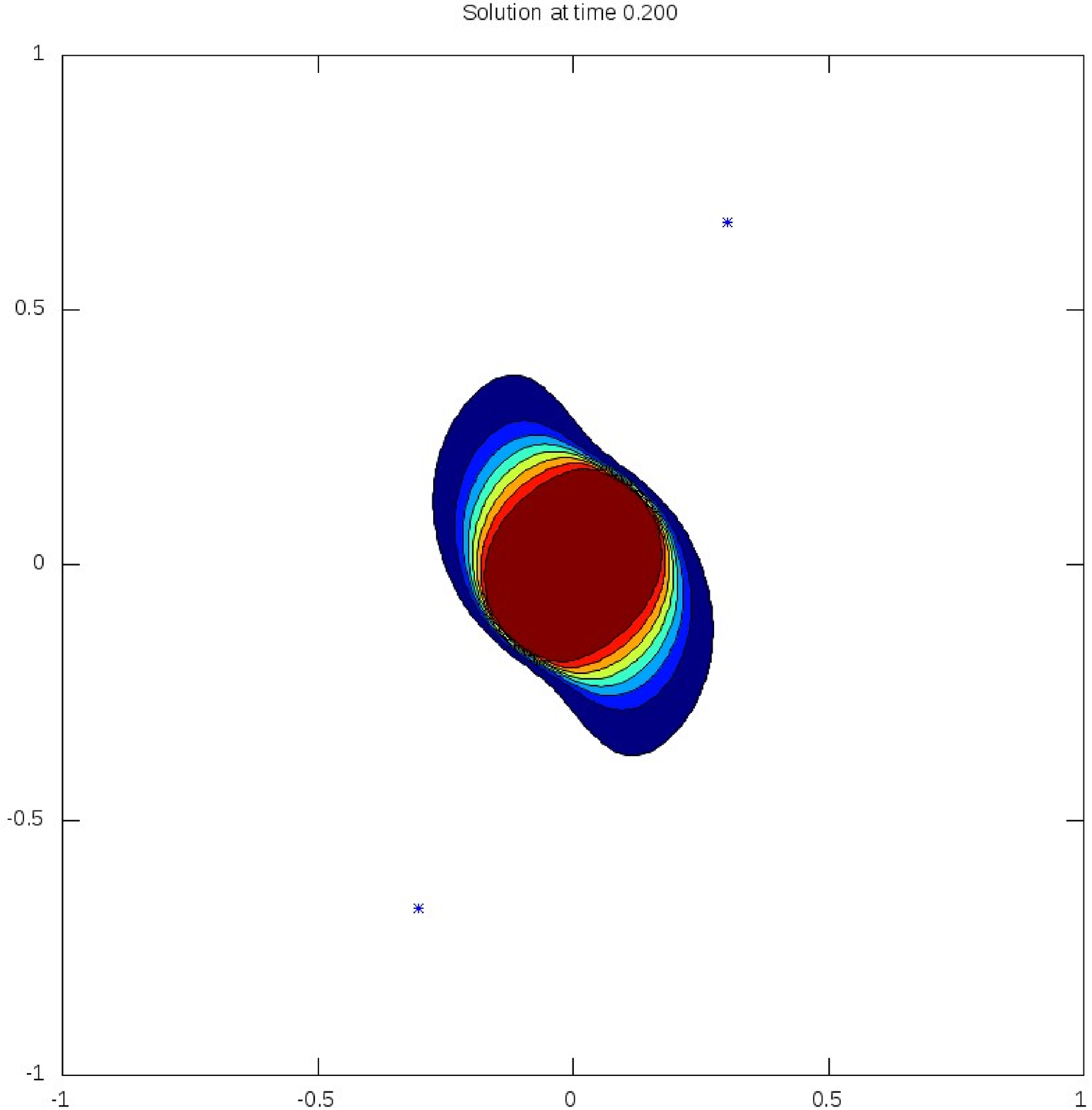}
  \Caption{Solution
    to~\eqref{eq:Problem}--\eqref{eq:dog_f}--\eqref{eq:ShepDog} at
    times $t = 0$, $t = 0.044$, $t = 0.067$, $t = 0.111$, $t = 0.156$,
    $t = 0.200$. Sheeps are initially uniformly distributed at the
    maximal density $R=1$ in the circumference centered at $(0,0)$
    with radius $0.2$. Dogs start moving from $(0.7, \, 0)$ and
    $(-0.7, \, 0)$, they succeed in confining the dispersion of the
    sheeps, at least for the tie interval considered.}
  \label{fig:Dog}
\end{figure}
Graphs of the corresponding solution are in Figure~\ref{fig:Dog}.

Merely technical modifications may allow to pass to various other
problems. For instance, dogs may be asked to constrain the movement of
all sheeps towards a certain area.

\subsection{Predator and Preys}
\label{subs:Hawk}

We consider here a predator attacking a group of preys. We think for
example at a hawk pursuing a flock of smaller birds or at a shark
attacking a group of sardines. Here, $\rho$ is the density of the
preys with $x\in \reali^3$, $p$ is now the pair $(P, V) \in \reali^6$,
where $P \in \reali^3$ is the position of the predator, $V \in
\reali^3$ is its speed and we postulate below an equation for the
acceleration $\ddot{P} =\dot V$ of the predator. Indeed, the framework
in Theorem~\ref{thm:main} allows to consider also second, or higher,
order ordinary differential equations for the single agents. The
initial density of the preys is assumed to have a compact, connected
support. The aim of the predator is to divide this connected group
into two smaller (disconnected) groups. Hence, its acceleration is
directed along the gradient of the preys' density, say $\ddot P =
\alpha \rho(t) *_x \nabla \eta$ for a suitable $\alpha >0$. The preys
have a speed $V_{\max} (1-\rho/R) V_0$, for a fixed $V_0 \in
\reali^2$, and modify it trying to escape from the predator. The
resulting speed of the preys is thus
\begin{equation}
  \label{eq:vPreys}
  v(t, x, \rho, p)
  =
  V_{\max} (1-\rho/R)
  \left(
    V_0
    +
    B \, e^{-C\norma{x-p(t)}} \, \left(x-p(t)\right)
  \right)
\end{equation}
where $B,C$ are positive constants. The former one is related to the
speed at which preys escape the predator and the latter to the
distance at which preys feel the presence of the predator.  Once
again, $v$ is maximal at zero density and vanishes at the maximal
density $R$, which means that the preys can not move when their
density is maximal.

\begin{lemma}
  \label{lem:preys}
  Let $N_x = 3$, $N_p = 6$, $N_r = 3$ and fix a positive $R$. Assume
  $v$ is as in~\eqref{eq:vPreys}, $\eta \in \Cc2(\reali^2, \reali)$
  with $\int_{\reali^2} \eta \d{x} =1$.  Denote $p = (P,V)$ and define
  \begin{equation}
    \label{eq:preys}
    f(t,x,\rho,p) = \rho \, v(t, x, \rho, p),
    \qquad
    \phi\left(t,
      \left[
        \begin{array}{c}
          P\\V
        \end{array}
      \right]
      ,r\right) =\left[
      \begin{array}{c}
        V\\
        r
      \end{array}
    \right] ,
    \qquad
    A\rho = \rho *_x\nabla \eta \,.
  \end{equation}
  Then, this setting fits in the framework of
  Corollary~\ref{cor:compact} as soon as $\bar \rho$ vanishes outside
  a compact set.
\end{lemma}

\begin{figure}[htpb]
  \centering
  \includegraphics[width=0.33\textwidth]{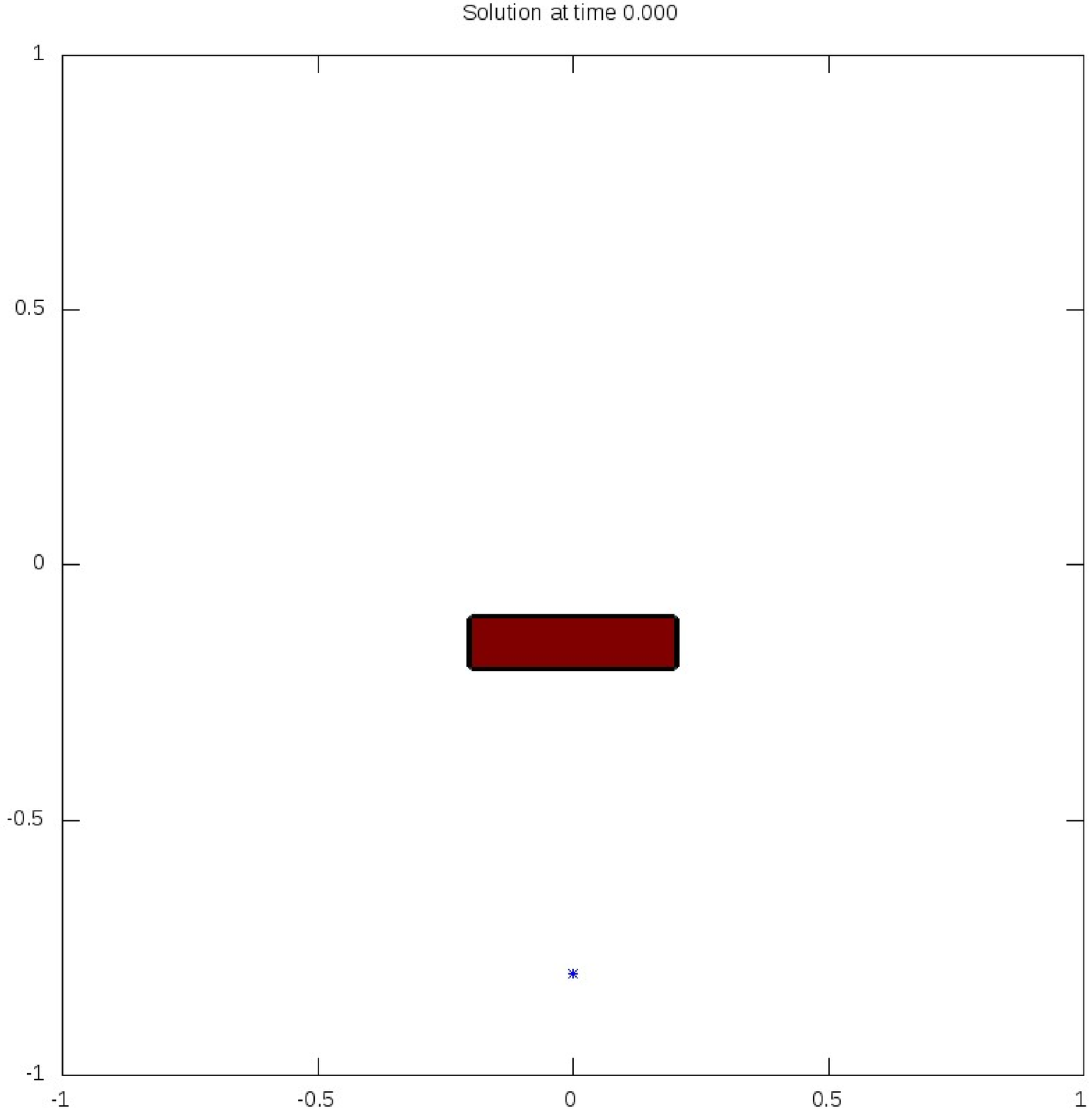}%
  \includegraphics[width=0.33\textwidth]{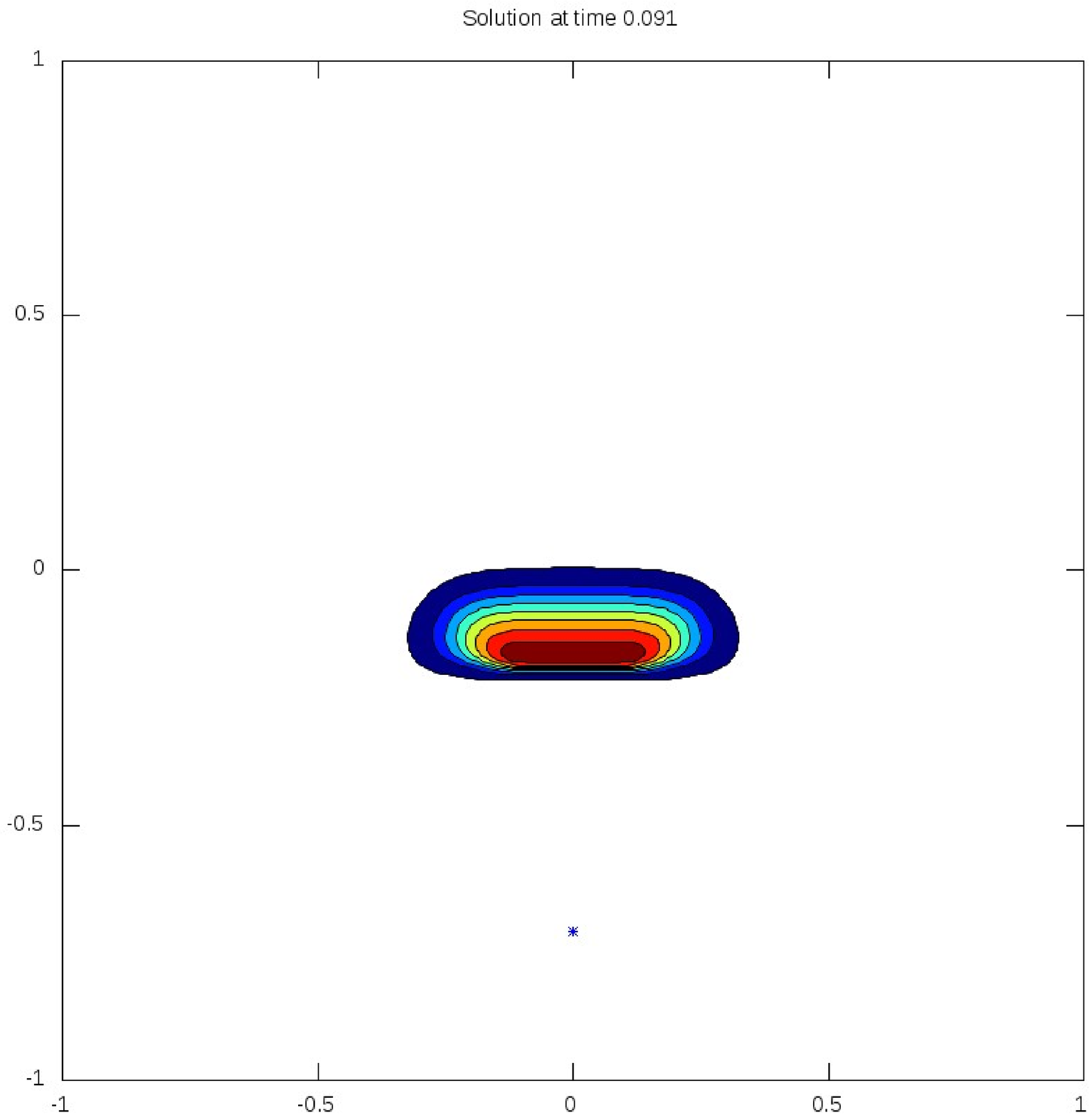}%
  \includegraphics[width=0.33\textwidth]{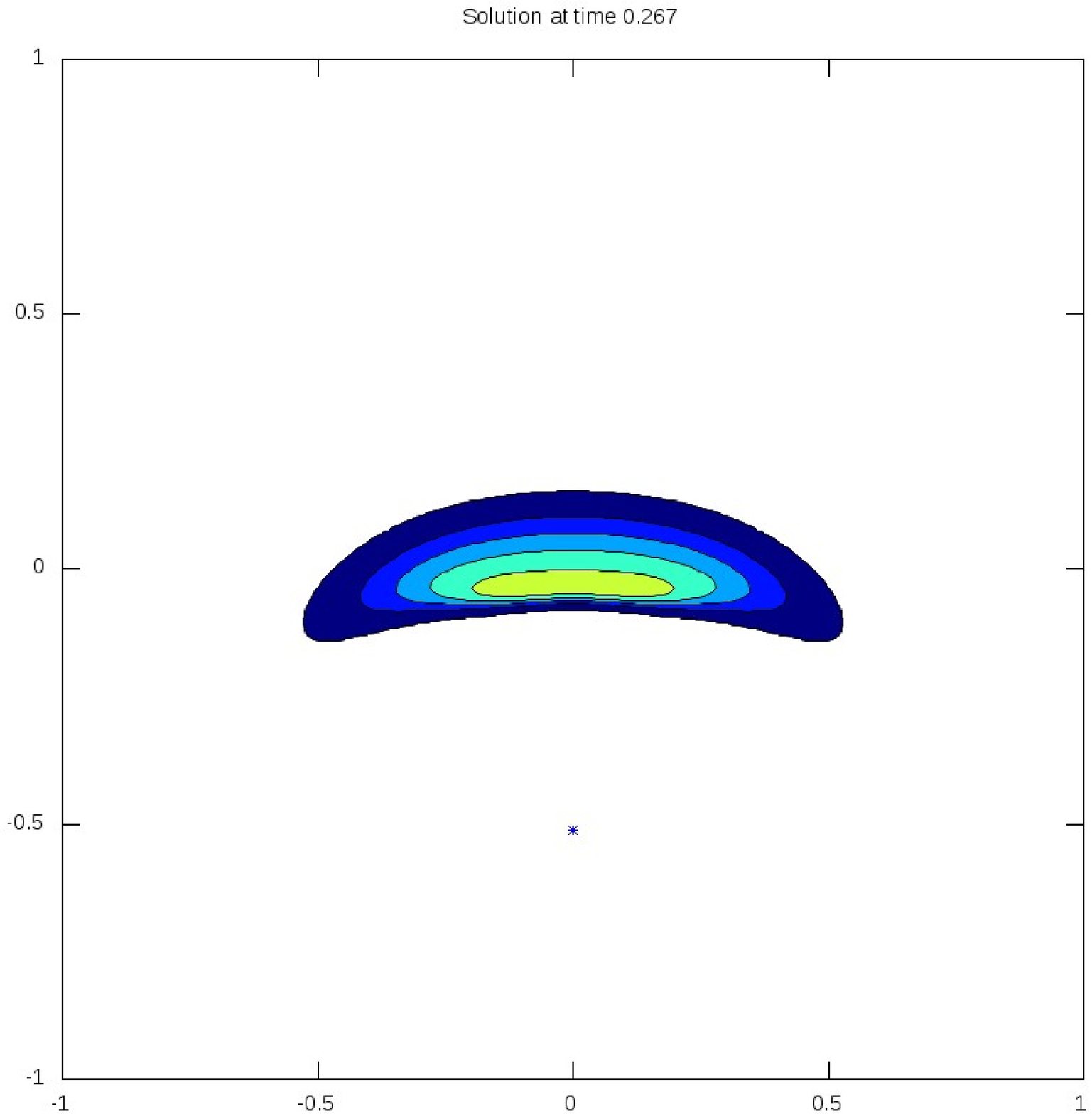}\\
  \includegraphics[width=0.33\textwidth]{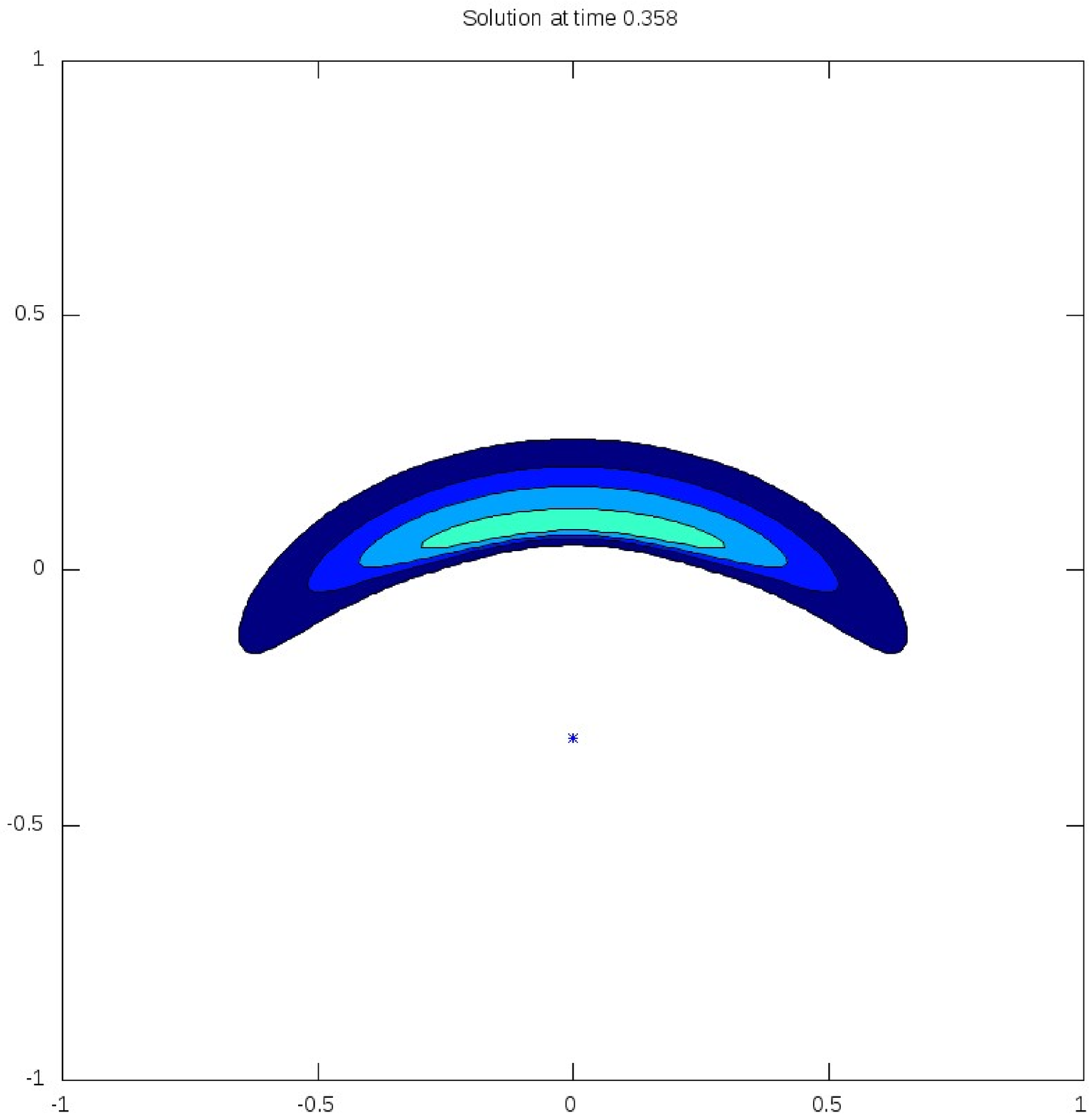}%
  \includegraphics[width=0.33\textwidth]{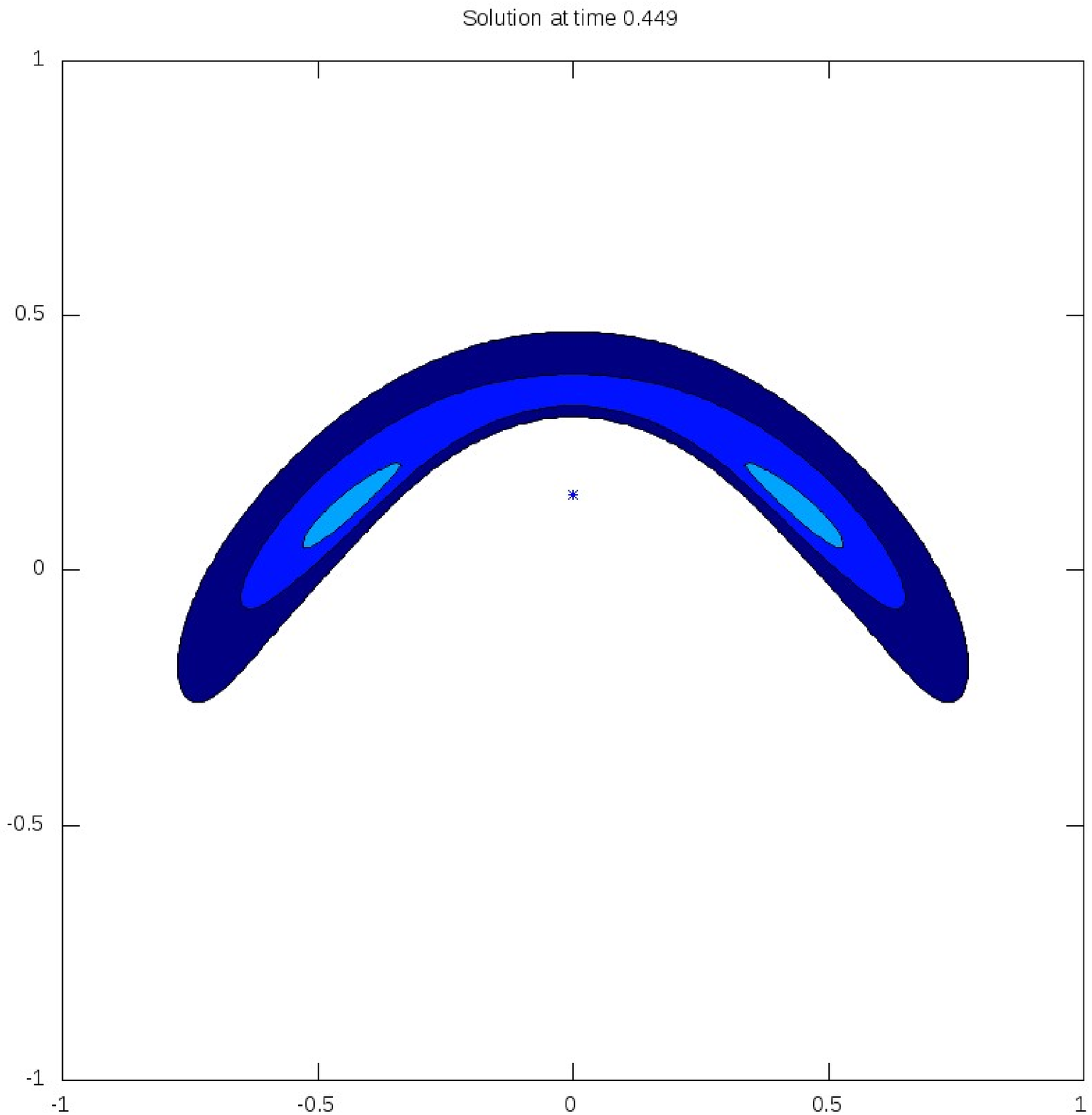}%
  \includegraphics[width=0.33\textwidth]{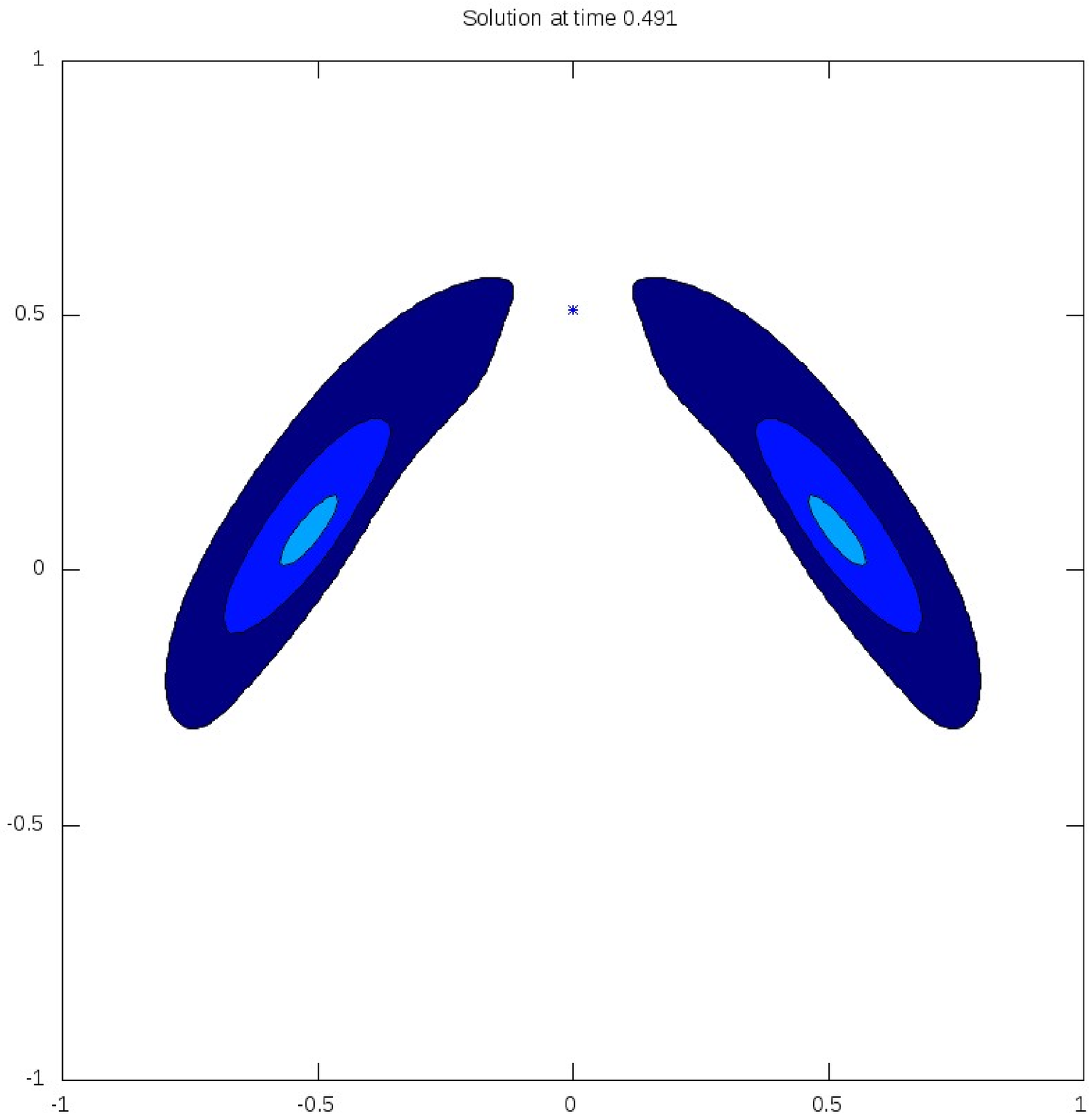}%
  \Caption{Solution obtained through the numerical integration
    of~\eqref{eq:Problem}--\eqref{eq:vPreys}--\eqref{eq:preys}--\eqref{eq:InitPrey}
    computed at times $0$, $0.091$, $0.267$, $0.358$, $0.449$ and
    $0.491$. Note that the predator succeeds in splitting the support
    of the preys.}
  \label{fig:Preys}
\end{figure}
\noindent\textbf{Numerical Integration:} For graphical purposes, we
limit the numerical integration to the 2D case. With reference
to~\eqref{eq:Problem}--\eqref{eq:vPreys}--\eqref{eq:preys}, we choose
the following parameters
\begin{equation}
  \label{eq:ParaPreys}
  \begin{array}{l}
    V_{\max} =  2
    \, , \quad
    C =  5.25
    \, , \quad
    V_0 =  [0 -0.5]^T
    \, , \quad
    B = 40
    \, , \quad
    A = 400
    \\
    \eta(x)
    =
    \displaystyle\frac{3}{\pi {r_p}^6}\,
    \left( \max \left\{0, {r_p}^2 - \norma{x}^2 \right\} \right)^2
    \qquad
    r_p = 0.5 \,.
  \end{array}
\end{equation}
and the initial datum
\begin{equation}
  \label{eq:InitPrey}
  P_o = \left[
    \begin{array}{c}
      0\\-0.8
    \end{array}
  \right]
  \,,\quad
  V_o =
  \left[
    \begin{array}{c}
      0 \\ 1
    \end{array}
  \right]
  \,,\quad
  \rho_o (x,y) = \chi_{\strut[-0.2, 0.2]}(x) \, \chi_{\strut [-0.2, -0.1]}(y) \,.
\end{equation}
The result is in Figure~\ref{fig:Preys}. Note that the predator
succeeds in splitting the support of the preys.

\section{Technical Details}
\label{sec:Tech}

Throughout, we let $W_{N_x} = \int_0^{\pi/2} (\cos \theta)^{N_x} \d
\theta$. We state below a Gr\"onwall--type lemma for later use.

\begin{lemma}
  \label{lem:Gronwall}
  Let the functions $\alpha \in \C0(I;\reali)$, $\beta \in
  \W{1}{1}(I;\reali)$, $\gamma \in \C0(I;\reali^+)$, $\Delta \in
  \C0(I;\reali)$ be such that
  \begin{displaymath}
    \Delta(t)
    \leq
    \alpha(t)
    \left(
      \beta(t) + \int_0^t \gamma(\tau) \, \Delta(\tau) \, \d\tau
    \right)\,.
  \end{displaymath}
  Then, for all $t \in I$,
  \begin{eqnarray*}
    \Delta(t)
    & \leq &
    \alpha(t)
    \left[
      \beta(0)
      \exp \left( \int_0^t \alpha(\tau)\, \gamma(\tau) \d\tau \right)
      +
      \int_0^t \beta'(\tau)
      \exp \left( \int_\tau^t \alpha(s) \, \gamma(s) \, \d{s} \right)
      \d\tau
    \right]
    \,.
  \end{eqnarray*}
\end{lemma}

\begin{proof}
  Using the following straightforward computations, we have:
  \begin{eqnarray*}
    \gamma(t)\Delta(t)
    &\leq &
    \alpha(t) \, \beta(t) \, \gamma(t)
    +
    \alpha(t) \, \gamma(t) \int_0^t \gamma(\tau) \Delta(\tau) \d{\tau}\,,
    \\
    \left(e^{-\int_0^t \alpha(\tau)\gamma(\tau)\d{\tau}}\int_0^t
      \gamma(\tau)\Delta(\tau)\d{\tau}\right)'
    &\leq &
    \alpha(t) \, \beta(t) \, \gamma(t) \,
    e^{-\int_0^t \alpha(\tau)\gamma(\tau)\d{\tau}}
    \,.
  \end{eqnarray*}
  Then, by integration we obtain
  \begin{eqnarray*}
    \int_0^t \gamma(\tau)\Delta(\tau)\d{\tau}
    &\leq &
    \int_0^t
    \alpha(t) \, \beta(t) \, \gamma(t) \,
    e^{\int_\tau^t \alpha(s)\gamma(s)\d{s}}
    \d\tau\,.
  \end{eqnarray*}
  Consequently, we have
  \begin{eqnarray*}
    \Delta(t)
    &\leq&
    \alpha(t)
    \left[
      \beta(t)
      +
      \int_0^t e^{\int_\tau^t \alpha(s)\gamma(s)\d{s}}
      \alpha(\tau) \beta(\tau) \gamma(\tau) \d{\tau}
    \right]\,.
  \end{eqnarray*}
  Integrating by part the last integral, we have finally the desired
  estimate.
\end{proof}

\begin{proofof}{Lemma~\ref{lem:pde}}
  This proof consists in applying to the scalar conservation law
  $ \partial_t \rho + \div_x f^*(t,x,\rho) =0$ with flux
  $f^*(t,x,\rho) = f\left(t, x, \rho, p(t)\right)$ first the classical
  Kru\v zkov result~\cite[Theorem~5]{Kruzkov} and then the stability
  estimates in~\cite{ColomboMercierRosini}.

  \smallskip

  To apply Kru\v zkov Theorem, it is sufficient to verify
  condition~\textbf{(H1)} in~\cite[Theorem~2.5]{ColomboMercierRosini}
  or the slightly weakened form in~\cite{MercierStability}. Note that:
  $f^*$ is $\C2$ in $x$ and $\rho$ by~\textbf{(f.\ref{it:f1})}, and is
  $\C0$ in $t$ by the regularity of $\pi$. This regularity is
  sufficient in the proof of~\cite[Theorem~2.5]{ColomboMercierRosini},
  see also~\cite[Remark~4 in \S~5]{Kruzkov}
  and~\cite{MercierStability}. Moreover, for any $t \in I$
  \begin{displaymath}
    \textbf{(f.\ref{it:f3})}
    \quad \Rightarrow \quad
    \partial_\rho f^*
    \in
    \L\infty ([0, t] \times \reali^{N_x} \times [0, R] ;
    \reali^{N_x})
    \mbox{ and }
    \div_x f^*
    \in
    \L\infty ([0, t] \times \reali^{N_x} \times [0, R] ;
    \reali)\, .
  \end{displaymath}
  Kru\v zkov Theorem can then be applied on any interval $[0, t]$.

  Observe that by~\textbf{(f.\ref{it:f2})}, the constant functions
  $\check\rho(t,x) \equiv 0$ and $\hat\rho(t,x) \equiv R$ solve the
  conservation law~(\ref{eq:HCL}), independently from $\pi$. Then, by
  the Maximum Principle~\cite[Theorem~3]{Kruzkov}, we have that any
  solution $\rho$ to~\eqref{eq:HCL} satisfies $\rho (t,x) \in [0, R]$
  for a.e.~$(t,x) \in I \times \reali^{N_x}$ and for all $\pi \in
  \C0(I; \reali^{N_p})$.

  To prove the $\L1$ continuity in time and the $\tv$ bound, we
  apply~\cite[Theorem~2.5]{ColomboMercierRosini} in the weaker form
  in~\cite{MercierStability}. To this aim, we verify
  also~\textbf{(H2)} therein on the time interval $[0, t]$, for any $t
  \in I$.  By~\textbf{(f.\ref{it:f5})} and the continuity of $\pi$,
  $\nabla_x \partial_\rho f^* \in \L\infty ([0, t] \times \reali^{N_x}
  \times [0, R]; \reali^{N_x\times N_x} )$.  Note also that,
  by~\textbf{(f.\ref{it:f6})}, $ \int_0^t \int_{\reali^{N_x}}
  \norma{\nabla_x \div_x f^*(\tau,x,\rho)}_{\L\infty} \d x \d \tau <
  +\infty$, with an upper bound that depends on $\pi$.

  \smallskip

  We denote below $\Omega_t=[0,t]\times\reali^{N_x}\times [0,R]\times
  K_t$ where $K_t$ is as above.
  By~\cite[Theorem~2.2]{MercierStability}
  or~\cite[Theorem~2.5]{ColomboMercierRosini} we obtain the estimate
  \begin{displaymath}
    \tv \left( \rho(t) \right)
    \leq
    \tv(\bar\rho) e^{\kappa_t t}
    +
    N_x W_{N_x}
    \int_0^t e^{\kappa_t (t-\tau)} \int_{\reali^{N_x}}
    \norma{\nabla_x \div_x f\left(\tau, x, \cdot, \pi(\tau)
      \right)}_{\L\infty([0,R])} \, \d x \, \d \tau
  \end{displaymath}
  where $\kappa_t = (2N_x+1) \norma{\nabla_x\partial_\rho
    f}_{\L\infty(\Omega_t)}$.  This implies~(\ref{eq:TV}).

  The $\L1$--continuity in time of $\rho$ follows
  from~\cite[Remark~2.4]{ColomboMercierRosini}, thanks
  to~\textbf{(f.\ref{it:f7})} and to the bound on the total variation,
  see also~\cite[Proof of Lemma~5.3]{ColomboHertyMercier}.

  To estimate the dependence of the solution from the initial datum,
  we check the hypotheses~\textbf{(H3)} in~\cite{MercierStability}
  or~\cite{ColomboMercierRosini} and
  apply~\cite[Theorem~2.3]{MercierStability}
  or~\cite[Theorem~2.6]{ColomboMercierRosini}.

  Let $f_1$ and $f_2$ satisfy~\textbf{(f.\ref{it:f1})}, $\ldots$,
  \textbf{(f.\ref{it:f6})}. Assume that $\pi_1$ and $\pi_2$ are in
  $\C0([0,t], \reali^{N_p})$. Let $f_1^*$ and $f_2^*$ be the
  corresponding compositions. With obvious notation, define $K = K^1_t
  \cup K^2_t$ and compute
  \begin{eqnarray*}
    & &
    \sup_{\tau\in [0,t],x\in \reali^{N_x}, \rho\in [0,R] }
    \modulo{
      \partial_\rho f_1^* \left(\tau, x, \rho, \pi_1(\tau) \right)
      -
      \partial_\rho f_2^* \left(\tau, x, \rho, \pi_2(\tau) \right)}
    \\
    &\leq &
    \norma{\partial_\rho f_1-\partial_\rho f_2}_{\L\infty(\Omega_t)}
    +
    \norma{\partial_\rho\nabla_p f_2}_{\L\infty(\Omega_t)} \;
    \norma{\pi_1-\pi_2}_{\L\infty([0,t])}
  \end{eqnarray*}
  which is bounded by~\textbf{(f.\ref{it:f3})}
  and~\textbf{(f.\ref{it:f8})}.

  To complete~\textbf{(H3)}, it remains only to estimate the quantity
  \begin{eqnarray*}
    & &
    \int_0^t \int_{\reali^{N_x}}
    \norma{
      \div_x \left(
        f_1\left(\tau,x, \cdot,\pi_1(\tau) \right)
        -
        f_2\left(\tau,x,\cdot,\pi_2(\tau) \right)
      \right)
    }_{\L\infty([0,R];\reali)} \d x \, \d \tau
    \\
    &\leq &
    \int_0^t \int_{\reali^{N_x}}
    \norma{\div_x(f_1-f_2) \left(\tau,
        x,\cdot,\pi_1(\tau) \right)}_{\L\infty([0,R])}\d{x}\d{\tau}
    \\
    & &
    +\int_0^t\int_{\reali^{N_x}} \norma{\nabla_p\div_x
      f_2(x)}_{\L\infty}\norma{\pi_1(\tau)-\pi_2(\tau)}\d{x}\d{\tau}
    \\
  \end{eqnarray*}
  which is bounded thanks to~\textbf{(f.\ref{it:f7})}
  and~\textbf{(f.\ref{it:f8})}. Now, we compare $\rho_1$ and $\rho_2$,
  obtaining
  \begin{align}
    \nonumber & \norma{(\rho_1-\rho_2)(t)}_{\L1}
    \\
    \nonumber \leq & \norma{\bar\rho_1-\bar\rho_2}_{\L1}
    \\
    & \nonumber + \!\left[ \frac{ e^{\kappa_t t}-1}{\kappa_t} \tv(\bar
      \rho ) + N_x W_{N_x} \!  \int_0^t \frac{ e^{\kappa_t
          (t-\tau)}-1}{\kappa_t} \int_{\reali^{N_x}}\!
      \norma{\nabla_x\div_x f_1\left(\tau, x,\cdot,
          \pi_1(\tau)\right)}_{\L\infty([0,R])}\! \d{x} \d{\tau}
    \right]
    \\
    \nonumber & \qquad \times \left( \norma{\partial_\rho
        f_1-\partial_\rho f_2}_{\L\infty(\Omega_t)} +
      \norma{\partial_\rho\nabla_p f_2}_{\L\infty(\Omega_t)}
      \norma{\pi_1-\pi_2}_{\L\infty([0,t])} \right)
    \\
    \nonumber & + \int_0^t \int_{\reali^{N_x}} \left( \norma{\div
        (f_1-f_2) \left(\tau, x,\cdot,\pi_1(\tau) \right)
      }_{\L\infty([0,R])} \right.
    \\
    \label{eq:explicit}
    & \left.  \qquad\qquad + \norma{\nabla_p \div_x f_2 (\tau,
        x,\cdot,\cdot)}_{\L\infty([0,R]\times K_t)}
      \norma{\pi_1(\tau)-\pi_2(\tau)}\right) \d{x}\d{\tau}
  \end{align}
  which gives the final estimate.
\end{proofof}

\begin{proofof}{Lemma~\ref{lem:edo}}
  By~\textbf{($\boldsymbol\phi$)}, we may apply~\cite[theorems~1
  and~2, Chap.~1]{Filippov} to~(\ref{eq:edo}) and get the local in
  time existence and uniqueness of the solution. The
  bound~(\ref{eq:SupBound}) follows from a standard application of
  Gr\"onwall Lemma and ensures that the solution can be extended to
  the whole interval $I$. Assume for simplicity that $\phi_1$ and
  $\phi_2$ satisfy~\textbf{($\boldsymbol{\phi})$} with the same
  function $C_\phi$.  Using the representation formula
  \begin{displaymath}
    p_i
    =
    \bar p_i
    +
    \int_0^t
    \phi_i\left( \tau, p_i(\tau), A\left(r_i(\tau)\right) (p_i(\tau)) \right)
    \d{\tau}\,,
  \end{displaymath}
  we get
  \begin{eqnarray*}
    & &
    \norma{(p_1-p_2)(t)}
    \\
    & \leq &
    \norma{\bar p_1-\bar p_2}
    +\!
    \int_0^t
    \norma{
      \phi_1 \left(
        \tau, p_1(\tau), A_1\left(r_1(\tau)\right) (p_1(\tau))
      \right)
      -
      \phi_2 \left(
        \tau, p_2(\tau), A_2\left(r_2(\tau)\right) (p_2(\tau))
      \right)
    } \!\d{\tau}
    \\
    & \leq &
    \norma{\bar p_1-\bar p_2}
    +
    \!\int_0^t
    \norma{(\phi_1-\phi_2)(\tau, p_1(\tau),
      A_1(r_1(\tau))(p_1(\tau)))} \d\tau
    \\
    & &
    \quad +
    \int_0^t
    C_\phi(\tau) \left(
      \norma{(p_1-p_2)(\tau)}
      +
      \norma{
        A_1\left(r_1(\tau)\right)(p_1(\tau))
        -
        A_2\left(r_2(\tau)\right)(p_2(\tau))}
    \right) \d{\tau}
    \\
    & \leq &
    \norma{\bar p_1-\bar p_2}
    +
    \int_0^t
    C_\phi(\tau)
    \left(1 + \norma{\nabla_p A_1(r_1)}_{\L\infty} \right)
    \norma{(p_1-p_2)(\tau)} \d{\tau}
    \\
    & &
    \quad
    +
    \int_0^t C_\phi (\tau)
    \left(
      \norma{A_1}_{\mathcal{L}(\L1,\W{1}{\infty})}
      \norma{(r_1-r_2)(\tau)}_{\L1}
      +
      \norma{A_1-A_2}_{\mathcal{L}(\L1,\W{1}{\infty})}
      \norma{r_2(\tau)}_{\L1}
    \right) \d\tau \\
    && +
    \int_0^t \norma{(\phi_1-\phi_2)(t, \cdot, \cdot)}_{\L\infty} d\tau\,.
  \end{eqnarray*}
  An application of Lemma~\ref{lem:Gronwall} with
  \begin{align*}
    \Delta(t) = & \norma{\bar p_1 - \bar p_2}\,,
    \\
    \alpha(t) = & 1\,,
    \\
    \beta(t) = & \norma{\bar p_1 - \bar p_2} + \int_0^t
    \norma{(\phi_1-\phi_2)(\tau,\cdot, \cdot)}_{\L\infty} \d\tau\,,
    \\
    \gamma(t) = & C_\phi(t) \left( 1 +
      \norma{A_1}_{\mathcal{L}(\L1,\W1\infty)} \norma{r_1}_{\L1}
    \right)
    \\
    & + \int_0^t C_\phi (\tau) \left[
      \norma{A_1}_{\mathcal{L}(\L1,\W{1}{\infty})}
      \norma{(r_1-r_2)(\tau)}_{\L1} +
      \norma{A_1-A_2}_{\mathcal{L}(\L1,\W{1}{\infty})}
      \norma{r_2(\tau)}_{\L1} \right] \!\d\tau .
  \end{align*}
  completes the proof of~(\ref{eq:est_edo}).
\end{proofof}

\begin{proofof}{Theorem~\ref{thm:main}}
  The proof is divided in several steps.

  \smallskip

  \noindent\textbf{1.~Local Existence.}
  Here we rely on an application of Banach Fixed Point Theorem. Fix
  first the initial data $\bar \rho \in (\L1 \cap \BV)(\reali^{N_x};
  [0,R])$ and $\bar p \in \reali^{N_p}$. Choose a positive $\hat T \in
  I$ and, motivated by~\eqref{eq:SupBound}, call
  \begin{displaymath}
    \delta
    =
    \left( \norma{\bar p} + 1 \right)
    e^{\int_0^{\hat T} C_\phi(\tau) \d\tau} -1\,.
  \end{displaymath}
  For any positive $\mathcal{R}$, with $\int_{\reali^{N_x}}\bar \rho
  \d{x} \leq \mathcal{R}$, and for any $T \in \left]0, \hat T\right]$,
  define the complete metric spaces and the distance
  \begin{eqnarray*}
    X_\rho
    & = &
    \left\{
      \rho \in \C0 \left([0,T]; \L1 (\reali^{N_x}; [0,R]) \right)
      \colon
      \sup_{t \in [0,T]}\int_{\reali^{N_x}} \rho(t,x) \d{x} \leq \mathcal{R}
    \right\}\,,
    \\
    X
    & = &
    X_\rho
    \times
    \C0 \left([0,T];B_{\reali^{N_p}}(0, \delta)\right)\,,
    \\
    d \left((\rho_1,p_1); (\rho_2,p_2) \right)
    & = &
    \sup_{t \in [0,T]}
    \norma{\rho_1(t) - \rho_2(t)}_{\L1}
    +
    \sup_{t \in [0,T]}
    \norma{p_1(t) - p_2(t)} \,.
  \end{eqnarray*}
  Define the map $\mathcal{T} \colon X \to X$ by $\mathcal{T}(r,\pi) =
  (\rho,p)$ if and only if $\rho$ and $p$ solve the problems
  \begin{equation}
    \label{eq:Contraction}
    \left\{
      \begin{array}{l}
        \partial_t \rho
        +
        \div_x f\left( t, x, \rho, \pi(t) \right)
        =
        0
        \\
        \rho(0,x) = \bar \rho(x)
      \end{array}
    \right.
    \quad \mbox{ and }\quad
    \left\{
      \begin{array}{l}
        \dot p
        =
        \phi\left( t, p, \left(Ar(t)\right) (p) \right)
        \\
        p(0) = \bar p \,.
      \end{array}
    \right.
  \end{equation}
  Note that both problems admit a unique solution, by
  lemmas~\ref{lem:pde} and~\ref{lem:edo}. Moreover, by the
  conservative form of the former problem in~\eqref{eq:Contraction},
  $\int_{\reali^{N_x}} \rho(t,x) \d{x} = \int_{\reali^{N_x}} \bar
  \rho(x) \d{x} \leq \mathcal{R}$, so that $\mathcal{T}$ is well
  defined. Moreover, Lemma~\ref{lem:edo} shows that the solution to
  the latter problem in~\eqref{eq:Contraction} is in $\W1\infty
  \left([0,T];B_{\reali^{N_p}}(0, \delta)\right) \subset \C0
  \left([0,T];B_{\reali^{N_p}}(0, \delta)\right)$.

  To prove that $\mathcal{T}$ is a contraction, fix $(r_1,\pi_1)$ and
  $(r_2,\pi_2)$ and call $(\rho_i, p_i) = \mathcal{T}(r_i,
  \pi_i)$. Then, define $K_{\hat T} = B_{\reali^{N_p}} (0,\delta)$ and
  apply Lemma~\ref{lem:pde} with $t=T$. Note that $K_T \subseteq
  K_{\hat T}$. The former problem in~\eqref{eq:Contraction} is then
  solvable in $\C0 \left([0,T]; \L1 (\reali^{N_x}; [0,R]) \right)$ and
  the stability estimate~(\ref{eq:dep_qup}) yields
  \begin{displaymath}
    \sup_{t \in [0,T]} \norma{\rho_1(t) - \rho_2(t)}_{\L1}
    \leq
    T \, \mathcal{C}(\hat T)
    \sup_{t \in [0,T]} \norma{\pi_1(t)-\pi_2(t)} \,.
  \end{displaymath}
  Apply now~(\ref{eq:est_edo})
  \begin{displaymath}
    \sup_{t \in [0,T]} \norma{p_1(t)-p_2(t)}
    \leq
    C_A
    \int_0^T C_\phi(\tau) \, e^{F(T) - F(\tau)} \d\tau
    \sup_{t \in [0,T]} \norma{r_1(t) - r_2(t)}_{\L1}\,,
  \end{displaymath}
  where $F$ is defined as in~(\ref{eq:F}) and can here be bounded as
  \begin{equation}
    \label{eq:FF}
    F(t)
    \leq
    \left( 1 + C_A \mathcal{R} \right) \int_0^t C_\phi(\tau)
    \d\tau \,.
  \end{equation}
  Hence,
  \begin{displaymath}
    d \left( \mathcal{T}(\rho_1, p_1), \mathcal{T}(\rho_2, p_2) \right)
    \leq
    \max\left\{
      T \, \mathcal{C}(\hat T) ,
      C_A (e^{F(T)}-1)
    \right\}
    \, d \left( (\rho_1,p_1), (\rho_2,p_2) \right) \,.
  \end{displaymath}
  Choose now a sufficiently small $T$ so that $\mathcal{T}$ is a
  contraction. Then, its unique fixed point is the unique solution
  to~(\ref{eq:Problem}) defined on the time interval $[0,T]$.

  \smallskip

  \noindent\textbf{2.~Global Uniqueness:}
  Let now $(\rho_1,p_1)$ and $(\rho_2,p_2)$ be two solutions to the
  same problem~(\ref{eq:Problem}) and defined at least on a common
  time interval $[0,\check T] \subseteq I$. Define
  \begin{displaymath}
    T^*
    =
    \sup \left\{
      T \in [0,\check T] \colon  (\rho_1,p_1)(t) = (\rho_2,p_2)(t)
      \mbox{  for all } t \in [0,T]
    \right\} \,.
  \end{displaymath}
  By the uniqueness of the fixed point, $(\rho_1,p_1)(t) =
  (\rho_2,p_2)(t)$ for all $t \in [0,T]$, so that the set in the right
  hand side above is not empty.  Repeat Step~1 with initial datum
  $(\bar\rho^*,\bar p^*) = (\rho_1,p_1)(T^*) = (\rho_2,p_2)(T^*)$,
  which is possible since $p$ is bounded on $[0,T^*]$ and
  $\tv(\bar\rho^*)$ is bounded, by~(\ref{eq:TV}). Thus, we obtain that
  $(\rho_1,p_1)(t) = (\rho_2,p_2)(t)$ also on a right neighborhood of
  $T^*$.  This contradicts the maximality of $T^*$, unless $T^* =
  \check T$.

  \smallskip

  \noindent\textbf{3.~Global Existence:}
  Define now
  \begin{displaymath}
    T_*
    =
    \sup \left\{
      T \in I \colon \exists
      \mbox{ a solution to~(\ref{eq:Problem}) defined on } [0,T]
    \right\}
  \end{displaymath}
  and assume that $T_* < +\infty$.  By~(\ref{eq:SupBound}), $p$ is
  bounded on $\left[0,T_*\right[$ and since
  \begin{displaymath}
    \norma{p(t_2) - p(t_1)}
    \leq
    \modulo{
      \int_{t_1}^{t_2}
      C_\phi(\tau) \left(  1 + \norma{p(\tau)} \right) d\tau
    }
    \leq
    \left( 1 + \sup_{t \in [0,T_*]} \norma{p(t)} \right)
    \modulo{\int_{t_1}^{t_2} C_\phi(\tau) \, d\tau} \,,
  \end{displaymath}
  $p$ is also uniformly continuous. Hence the limit $p_* = \lim_{t \to
    T_*^-} p(t)$ exists and is finite.

  Apply now Lemma~\ref{lem:pde} on the interval $[0,T_*]$, obtaining
  that the solution $\rho$ to~(\ref{eq:HCL}) is defined on all
  $[0,T_*]$ and, together with $p$, also
  solves~(\ref{eq:Problem}). Now, we repeat Step~1 with initial datum
  $(\bar\rho_*, \bar p_*) = (\bar\rho, \bar p)(T_*)$, which is
  possible thanks to~(\ref{eq:TV}). In turn, this allows to extend
  $(\bar\rho, \bar p)$ to a right neighborhood of $T_*$.  This
  contradicts the maximality of $T_*$, unless $T_* = T_{\max}$.

  \smallskip

  \noindent\textbf{4.~Stability Estimates:}
  Fix $t>0$ and let $\tau \in [0,t]$. Let $\mathcal{R} \geq \max
  \left\{ \int_{\reali^{N_x}} \bar \rho_1 \d{x}, \int_{\reali^{N_x}}
    \bar \rho_2 \d{x} \right\}$.  Then, by~(\ref{eq:dep_qup})
  and~(\ref{eq:est_edo}), the solutions to~\eqref{eq:Problem2} satisfy
  \begin{eqnarray*}
    & &
    \norma{(\rho_1-\rho_2)(t)}_{\L1}
    \\
    & \leq &
    \norma{\bar\rho_1 - \bar\rho_2}_{\L1}
    +
    t \, \mathcal{C}(t)
    \Bigl[
    \norma{p_1-p_2}_{\L\infty([0,t])} +
    \norma{\partial_\rho(f_1-f_2)}_{\L\infty(\Omega_t)}
    \\
    & &
    \qquad\qquad
    +
    \norma{\div(f_1 - f_2)}_{\L1(\reali^{N_x}) \times \L\infty([0,t]
      \times [0,R] \times K_t)}
    \Bigr] \,,
    \\
    & &
    \norma{(p_1-p_2)(t)}
    \\
    & \leq &
    \displaystyle
    e^{F(t)} \norma{\bar p_1 - \bar p_2}
    +
    \int_0^t e^{F(t)-F(\tau)}
    \norma{(\phi_1 - \phi_2)(\tau, \cdot, \cdot)}_{\L\infty} \d\tau
    \\
    & &
    +
    \displaystyle
    \int_0^t \! e^{F(t)-F(\tau)} C_\phi(\tau)
    \left(
      C_A
      \norma{(\rho_1-\rho_2)(\tau)}_{\L1}
      +
      \mathcal{R}
      \norma{A_1-A_2}_{\mathcal{L}(\L1,\W{1}{\infty})}
    \right)
    \d\tau .
  \end{eqnarray*}
  with $\mathcal{C}$ as in Lemma~\ref{lem:pde}, $F$ as
  in~(\ref{eq:FF}), $K_t=B(0, \delta_t)$ and $\delta_t = \left(
    \norma{\bar p} + 1 \right) e^{\int_0^t C_\phi(\tau) \d\tau} -1$.
  Insert now the former estimate in the latter one and apply
  Lemma~\ref{lem:Gronwall} with
  \begin{eqnarray*}
    \Delta
    & = &
    \norma{(p_1 - p_2)(t)}\,,
    \\
    \alpha(t)
    & = &
    e^{F(t)}\,,
    \\
    \beta(t)
    & = &
    \norma{\bar p_1 - \bar p_2}
    +
    \frac{\mathcal{R}}{1+C_A\mathcal{R}} \left( 1 - e^{-F(t)} \right)
    \norma{A_1-A_2}_{\mathcal{L}(\L1,\W1\infty)}
    \\
    & &
    +
    \int_0^t \norma{(\phi_1 - \phi_2)(\tau, \cdot, \cdot)}_{\L\infty}
    e^{-F(\tau)} \d\tau
    +
    C_A \int_0^t e^{-F(\tau)} C_\phi(\tau) \norma{\bar \rho_1 - \bar
      \rho_2}_{\L1} \d\tau
    \\
    & &
    +
    C_A
    \int_0^t \tau \, \mathcal{C}(\tau) \, C_\phi(\tau) \, e^{-F(\tau)}
    \\
    & &
    \qquad
    \times
    \left(
      \norma{\partial_\rho(f_1-f_2)}_{\L\infty(\Omega_\tau)}
      +
      \norma{\div(f_1 - f_2)}_{\L1(\reali^{N_x}) \times \L\infty([0,\tau]
        \times [0,R] \times K_\tau)}
    \right)\,,
    \\
    \gamma(t)
    & = &
    C_A \, t \, C_\phi(t) \, \mathcal{C}(t)
    e^{F(t)}\,,
  \end{eqnarray*}
  obtaining, with $\mathcal{H}(\tau,t) = \exp \int_\tau^t C_\phi(s)
  \left( 1 + C_A\, \mathcal{R} + C_A \, s \, \mathcal{C}(s) \right)
  \d{s}$,
  \begin{eqnarray*}
    \norma{p_1-p_2}
    & \leq &
    \left(
      \exp \left(
        F(t)
        +
        C_A
        \int_0^t \tau C_\phi(\tau) \mathcal{C}(\tau) \d\tau
      \right)
    \right)
    \norma{\bar p_1 - \bar p_2}
    \\
    & &
    +
    \left( \int_0^t \mathcal{H}(\tau,t) \d\tau \right)
    \norma{\phi_1-\phi_2}_{\L\infty([0,t]\times K_t \times [0,C_A]}
    \\
    & &
    +
    \left(
      \mathcal{R}  \int_0^t C_\phi(\tau) \mathcal{H}(\tau,t) \d\tau
    \right)
    \norma{A_1-A_2}_{\mathcal{L}(\L1,\W1\infty)}
    \\
    & &
    +
    \left(
      C_A \int_0^t C_\phi(\tau) \, \mathcal{H}(\tau,t) \d\tau
    \right)
    \norma{\bar\rho_1-\bar\rho_2}_{\L1}
    \\
    & &
    +
    \left(
      C_A \int_0^t \tau C_\phi(\tau) \mathcal{C}(\tau)
      \mathcal{H}(\tau,t) \d\tau
    \right)
    \\
    & & \quad
    \times
    \left[
      \norma{\partial_\rho(f_1-f_2)}_{\L\infty([0,R]\times\reali^{N_x}\times
        K_t)}
      +
      \norma{\div_x(f_1-f_2)}_{\L1(\reali^{N_x})\times
        \L\infty([0,R]\times K_t)}
    \right].
  \end{eqnarray*}
  Then, we immediately get the other bound
  \begin{eqnarray*}
    &&\norma{\rho_1-\rho_2}_{\L1}\\
    & \leq &
    \norma{\bar \rho_1 - \bar \rho_2}
    \left(
      1 + t \mathcal{C}(t)
      \exp \left(
        F(t)
        +
        C_A
        \int_0^t \tau C_\phi(\tau) \mathcal{C}(\tau) \d\tau
      \right)
    \right)
    \\
    & &
    +
    \left(
      t \mathcal{C}(t)
      \int_0^t \mathcal{H}(\tau,t) \d\tau \right)
    \norma{\phi_1-\phi_2}_{\L\infty([0,t]\times K_t \times [0,C_A]}
    \\
    & &
    +
    \left(
      \mathcal{R}  t \mathcal{C}(t)
      \int_0^t C_\phi(\tau) \mathcal{H}(\tau,t) \d\tau
    \right)
    \norma{A_1-A_2}_{\mathcal{L}(\L1,\W1\infty)}
    \\
    & &
    +
    C_A t \mathcal{C}(t)
    \exp \left(
      F(t)
      +
      C_A \left( 1+\tv(\bar\rho_1)\right)
      \int_0^t \tau C_\phi(\tau) \mathcal{C}(\tau) \d\tau
    \right)
    \norma{\bar p_1-\bar p_2}_{\L1}
    \\
    & &
    +
    t \mathcal{C}(t)
    \left(
      1
      +
      C_A \int_0^t \tau C_\phi(\tau) \mathcal{C}(\tau)
      \mathcal{H}(\tau,t) \d\tau
    \right)
    \\
    & & \quad
    \times
    \left(
      \norma{\partial_\rho(f_1-f_2)}_{\L\infty([0,R]\times\reali^{N_x}\times
        K_t)}
      +
      \norma{\div_x(f_1-f_2)}_{\L1(\reali^{N_x})\times
        \L\infty([0,R]\times K_t)}
    \right)
  \end{eqnarray*}
  completing the proof.
\end{proofof}

We need below the following consequence of Kru\v zkov
Theorem~\cite[Theorem~5]{Kruzkov}.

\begin{proposition}
  \label{prop:FiniteSpeed}
  Let $N_x \in \naturali$ and $T>0$. Consider the conservation law
  \begin{equation}
    \label{eq:hcl}
    \left\{
      \begin{array}{l}
        \partial_t \rho + \div_x \bar f(t,x,\rho) =0
        \\
        \rho(t,0) = \bar \rho
      \end{array}
    \right.
  \end{equation}
  with $\bar f \in \C0([0,T] \times \reali^{N_x} \times \reali;
  \reali^{N_x})$; $\partial_\rho \bar f$, $\partial_\rho \nabla_x \bar
  f$ and $\nabla_x^2 \bar f$ continuous wherever defined;,
  $\partial_\rho \bar f, \, \div_x \bar f \in \L\infty([0,T] \times
  \reali^{N_x} \times [-H, H])$ for all $H>0$. Assume that $\bar\rho
  \in (\L1 \cap \L\infty)(\reali^{N_x};\reali)$ is such that $\bar
  \rho (x) = 0$ for a.e.~$x \in \reali^{N_x} \setminus
  B_{\reali^{N_x}}(0, d)$ for a given $d>0$. Moreover, $\bar f(t, x,
  0) = 0$ for all $t \in [0,T]$ and $x \in \reali^{N_x}$.  Call $\rho$
  the Kru\v zkov solution to~\eqref{eq:hcl} and let $K = \sup_{t \in
    [0,T]} \norma{\rho(t)}_{\L\infty(\reali^{N_x})}$. Then, for all $t
  \in [0,T]$, $\rho(t,x) =0$ for a.e.~$x \in \reali^{N_x} \setminus
  B_{\reali^{N_x}}(0, d+ Vt)$, where $V = \norma{\partial_\rho \bar
    f}_{\L\infty([0,T] \times \reali^{N_x} \times [-K,K])}$.
\end{proposition}

Above, $\bar f$ is assumed to satisfy the usual Kru\v zkov conditions,
see~\cite[\textbf{(H1)}]{MercierStability},
or~\cite{ColomboMercierRosini, Kruzkov}. The proof essentially relies
on~\cite[Theorem~1]{Kruzkov}.

\begin{proofof}{Proposition~\ref{prop:FiniteSpeed}}
  Choose an $x \in \reali^{N_x} \setminus
  B_{\reali^{N_x}}(0,d+Vt)$. Let $\delta>0$ be such that
  $B_{\reali^{N_x}}(x,\delta) \cap B_{\reali^{N_x}}(0,d+Vt) =
  \emptyset$, so that $B_{\reali^{N_x}}(x,\delta+Vt) \cap
  B_{\reali^{N_x}}(0,d) =
  \emptyset$. Applying~\cite[Theorem~1]{Kruzkov}, with $u=\rho$ and
  $v=0$, we have that
  \begin{displaymath}
    \int_{B_{\reali^{N_x}}(x,\delta)} \modulo{\rho(t,x)} \d{x}
    \leq
    \int_{B_{\reali^{N_x}}(x,\delta+Vt)} \modulo{\bar \rho(x)} \d{x}
    =
    0
  \end{displaymath}
  hence $\rho(t)$ vanishes a.e.~outside $B_{\reali^{N_x}}(0,d+Vt)$.
\end{proofof}

\begin{proofof}{Corollary~\ref{cor:compact}}
  Fix any positive $T \in I$. Let $d$ be such that $\bar\rho$ vanishes
  outside $B_{\reali^{N_x}}(0,d)$ and call $\mathcal{K} =
  B_{\reali^{N_x}}(0, d+VT)$. Let $\chi \in \Cc\infty(\reali, [0,1])$
  be such that $\chi(x) = 1$ for all $x \in \mathcal{K}$. Define the
  convolution in the space variable $f_* = f*_x \chi$, so that $f^*$
  has compact support in $x$. Then, thanks also to the \emph{a priori}
  bound~\eqref{eq:SupBound}, $f^*$ satisfies~\textbf{(f)} on the
  interval $[0,T]$. Hence to the problem
  \begin{displaymath}
    \left\{
      \begin{array}{l}
        \partial_t \rho
        +
        \div_x f^* \left(t, x, \rho, p(t)\right)
        =
        0
        \\
        \dot p
        =
        \phi\left( t, p, \left(A \rho(t)\right) (p) \right)
        \\
        \rho(0,x) = \bar \rho(x)
        \\
        p(0) = \bar p
      \end{array}
    \right.
  \end{displaymath}
  Theorem~\ref{thm:main} can be applied, yielding the existence and
  uniqueness of a solution $(\rho,p)$ in the sense of
  Definition~\ref{def:sol} defined on all the interval $[0,T]$. Let
  now $\bar f(t, x, \rho) = f^*\left( t, x, \rho, p(t) \right)$. Then,
  $\rho$ is a Kru\v zkov solution to~\eqref{eq:hcl} and by
  Proposition~\ref{prop:FiniteSpeed} its support is contained in
  $\mathcal{K}$, for all $t \in [0,T]$. Thereof ore, on the same time
  interval, by the definition of $f^*$, $(\rho,p)$ is the unique
  solution also to~\eqref{eq:Problem}, always according to
  Definition~\ref{def:sol}. The rest of the proof easily follows.
\end{proofof}

{\small

  \bibliographystyle{abbrv}

  \bibliography{luminy4}

}

\end{document}